\documentclass[12pt]{elsarticle}

\usepackage{amssymb}
\usepackage[dvipsnames]{xcolor}
\usepackage{array}
\usepackage{amsmath}
\usepackage{subfigure}
\usepackage{graphicx}
\usepackage{float}
\usepackage{setspace}
\usepackage{booktabs}
\usepackage{lscape}
\usepackage{algorithm} 
\usepackage{algpseudocode}
\usepackage{rotating}
\usepackage{nomencl}

\journal{arXiv}

\begin{document}

\begin{frontmatter}

\title{An areal continuum model for mixed traffic }

\author[inst1]{Nandan Maiti}

\affiliation[inst1]{organization={Univ. Eiffel, ENTPE, LICIT-ECO7},
            addressline={Lyon}, 
            postcode={69500},
            country={France}}

\author[inst2]{Bhargava Rama Chilukuri}
\affiliation[inst2]{organization={Department of Civil Engineering, Indian Institute of Technology},
            addressline={ Madras}, 
            city={Chennai},
            postcode={600036}, 
            state={Tamil Nadu},
            country={India}}

\begin{abstract}
A novel continuum model has been developed to address the vehicle size heterogeneity in mixed traffic. By incorporating the principle of vehicle area conservation, a new set of traffic flow variables centered on the concept of vehicle area has been introduced. These variables, namely areal flow and areal density, exhibit the remarkable characteristic of being conserved across time and space and seamlessly related to the traditional space mean speed. Next, the study successfully established a bi-variate relationship among these newly introduced area-based traffic variables based on empirical data from three locations for stream and mode-specific. It is shown that the conventional solution methods for the hyperbolic conservation laws remain applicable to the proposed continuum model. Finally, a multi-class cell transmission model numerical scheme is developed and used to illustrate the performance of the proposed model in replicating the seepage and platoon dispersion behavior in mixed traffic conditions.
\end{abstract}




\end{frontmatter}

\section{Introduction}
 In recent years, there has been an increasing interest in studying traffic flow characteristics for mixed traffic conditions. One of the key considerations in mixed traffic flow modeling is the heterogeneity of vehicle types and lane-free movement. Mixed traffic conditions are generally characterized by a mix of vehicles such as cars, motorcycles, bicycles, trucks, buses, light commercial vehicles, and other informal modes of transportation \citep{dhamaniya2013concept,maiti2023empirical}. These diverse vehicle types exhibit different behavior patterns, speeds, and flow characteristics, requiring modeling approaches to capture this heterogeneity \citep{mallikarjuna2011heterogeneous,maiti2023does,maiti2023empirical,maiti2023estimation}.
\par Fundamental diagrams (FDs) are one of the critical models for traffic flow analysis.  Empirical plots of heterogeneous traffic often exhibit significant scatterings despite carefully identifying the steady states \citep{maiti2023empirical}. While nonstationary behaviors and other sources of randomness like acceleration capability, maneuverability, and dynamic influence zone of aggressive vehicles contribute to this scatter, the interaction between different vehicle classes is a crucial factor that requires careful investigation. Traditionally, the impact of each vehicle class on traffic state variables is quantified by converting all classes to a reference class using Passenger Car Units (PCUs). 
\par The concept of PCU is widely accepted to standardize heterogeneous traffic into a homogeneous equivalent. HCM 2000 defined PCU for lane-based conditions, and several researchers \citep{chandra2003effect,krishnamurthy2008effect,biswas2017estimation} embraced this concept and expanded their analysis to determine the PCU values of various vehicle categories, including heavy vehicles, in the context of heterogeneous traffic conditions. The most popular speed and headway-based PCU estimation method of \citet{chandra2003effect} is compatible with heterogeneous traffic conditions. Numerous studies and capacity manuals have suggested static PCU values \citep{fan1990passenger,adnan2014passenger}, while more recent research has highlighted the dynamic nature of PCUs \citep{arasan2010simulation,dhamaniya2016conceptual,biswas2017estimation}. These studies have revealed that the PCU of a vehicle can vary based on the surrounding traffic and other factors such as vehicle speeds. Nonetheless, there are recognized drawbacks in utilizing PCU to represent heterogeneous traffic flows, including challenges related to model parameter calibration \citep{lee1998origin,wang2010representing,dhamaniya2013concept}, class-specific speed estimations \citep{Benzoni-Gavage02ann-populations,Ngoduy2007}, neglecting mixed vehicle interactions \citep{bouadi2016effect}  and the absence of creeping modeling \citep{VanLint2008,Ngoduy2007}. Therefore, accurately defining state-dependent PCU values that reflect the dynamics among vehicle classes with vastly different characteristics remains challenging.
\par  The inter-vehicle heterogeneity discussed in \citet{ossen2006multi} refers to the variances in behavior or capabilities between different classes of vehicles and drivers against PCUs conversion. Instead of converting heterogeneity to equivalent homogeneous traffic, researchers proposed multi-class models considering the heterogeneity of vehicles and drivers by dividing them into classes with distinct properties \citep{daganzo2002behavioral1,WONG2002827,gupta2007new,logghe2008multi}.  Multi-class models are developed based on the distinguishing class-specific speed  $v^i(.)$ in different regimes. Multi-class models estimate the individual vehicle class velocity by considering the aggregate and vehicle-wise variables. The effective density, in turn, is estimated by considering the principles of static \citep{Ngoduy2007,WONG2002827,zhu2003numerical} and dynamic PCU conservation \citep{Benzoni-Gavage02ann-populations, van2008fastlane,Nair2011,van2014new,fan2015heterogeneous}. Consequently, the limitations associated with the PCU conversion align with estimating multiclass fundamental diagrams. Also, incorporating more than two vehicle classes introduces increased complexity in defining regimes, identifying boundaries, and analyzing their properties due to the mixing of class-wise regimes and the corresponding states.  
\par Traditionally, mixed traffic conditions are characterized based on vehicle counts through PCU conversion or occupancy as indicators of traffic density. However, these traditional approaches fail to capture the complex nature of mixed traffic scenarios, where vehicles of different sizes and types coexist. \citet{khan1999modeling} argued that area measures address the limitation of the traditional variable by considering the spatial occupancy of vehicles, incorporating the physical area they occupy on the roadway. Later, \citet{mallikarjuna2006area} defined area occupancy, representing the temporal proportion during which the vehicle's area occupies the detection zone area. \citet{arasan2008measuring} validated the area occupancy concept to measure traffic concentration in heterogeneous traffic using a simulation technique. \citet{mohan2021multi} modified this concept by introducing perceived area occupancy, which incorporates the spatial arrangement and percentage composition of vehicles. In a recent study, \citet{suvin2018modified} proposed a modified approach that introduces area density and flow for mixed traffic. They achieved this by considering the width of vehicles and extending Edie's definition \citep{edie1963discussion} to incorporate the lateral dimension with the time-space region. However, representing heterogeneous traffic by considering only vehicle width \citep{suvin2018modified} or converting to equivalent density using vehicle length \citep{van2008fastlane,fan2015heterogeneous} instead of the whole vehicle area can lead to problems in capacity estimation, traffic flow characteristics, lane utilization, lane changing behavior, and congestion modeling in mixed traffic conditions.
\par The dynamic traffic flow model integrates fundamental diagrams and the continuum model to describe traffic behavior. The classical macroscopic traffic model, proposed by \citet{lighthill1955kinematic,richards1956shock}, establishes a dynamic equation for density in the time-space domain. Consequently, most models for heterogeneous traffic incorporate class-specific continuum models and multi-class fundamental diagrams \citep{WONG2002827,gupta2007new,logghe2008multi,fan2015heterogeneous} or stream conservation equations utilizing PCU-based fundamental diagrams \citep{krishnamurthy2008effect,biswas2017estimation,chandra2003effect,arasan2010simulation,biswas2017estimation,dhamaniya2016conceptual}. However, the headway or speed-based PCU value change during every time step, changing the effective density defying the concept of vehicle conservation. In some compositions, the effective density may even exceed jam density with the new PCU values. As the definition of the effective density implies, it cannot be conserved over time-space; hence, there is no equivalent of the continuity equation \citep{van2014new}. Additionally, due to the frequently changing compositions, the PCU conversion for mixed vehicles is computationally expensive, inaccurate, and difficult to estimate. Several continuum models used density alternatives in mixed traffic, such as those proposed by \citet{Nair2011, qian2017modeling, mohan2017heterogeneous}, have adopted similar concepts. However, most of the literature replaces the density with the alternative density variables instead of demonstrating the conservation form of those variables in the time-space domain. 
Thus, the literature lacks an appropriate dynamic modeling framework for mixed traffic conditions. Additionally, there's a notable absence of adequate macroscopic variables and continuum models beyond area occupancy for mixed traffic, rendering multi-class models potentially inadequate. Addressing these gaps is crucial for a nuanced comprehension of mixed traffic dynamics.
  Therefore, a continuum model is needed to accurately conserve parameters like effective density over time-space to address these challenges in Eulerian formulations. In other words, since the principle of vehicle number conservation may not be suitable for mixed traffic due to the vehicle class heterogeneity, a continuity equation based on vehicle area conservation is a potential direction, but has not been explored in the literature. 
\par The need for a continuum model that conserves vehicle area rather than vehicle number arises from the heterogeneous nature of traffic, where vehicles vary in size and shape. Traditional traffic flow models fail to capture the spatial occupancy of vehicles and the impact of different vehicle sizes on traffic dynamics. Thus an area-based continuum model would consider the physical space occupied by vehicles, allowing for a more accurate characterization of the mixed-traffic environment.

\subsection{Need for the study and major contributions}
The research gap in the given topic lies in the need for a more comprehensive and accurate modeling approach for mixed traffic conditions Specifically, there is a need to address the challenges associated with modeling mixed traffic in terms of vehicle types, sizes, and behaviors. The existing methods of converting heterogeneous traffic into equivalent homogeneous traffic using PCUs have several drawbacks, including challenges in parameter calibration, neglecting mixed vehicle interactions, and the dynamic nature of PCU values. First, suitable traffic flow variables that appropriately characterize the traffic conditions accounting for the physical area that vehicles occupy on the roadway and then a continuum model that accurately conserves effective density is necessary to improve the modeling and understanding of mixed traffic flow dynamics. 

The major contributions of this paper are:
\begin{enumerate}

    \item This paper presents new traffic flow variables and their definition based on vehicle area: areal flow and areal density. Furthermore, it highlights the distinctive characteristics of these metrics, demonstrating their compatibility with conventional traffic flow variables and theories.
    \item This research validated steady-state relationships of the proposed variables using existing fundamental diagram models utilizing empirical data sets.
    \item This paper developed an areal continuum model for mixed traffic based on the principle of vehicle area conservation.    \item  This research introduces a multi-class cell transmission model to demonstrate its capability to replicate mixed traffic phenomena such as overtaking, platoon dispersion, etc.      
\end{enumerate}
\par The rest of the paper is organized as follows: Section II describes the derivation of the continuum model using the principle of area conservation and proposed area-based traffic flow variables. Section III analyzes existing stream-based fundamental relationships and explores multi-class fundamental diagrams for mixed traffic incorporating new variables. Section IV discusses the applicability of existing numerical schemes to the proposed continuum model and proposes a multiclass cell transmission model to explain platoon dispersion and overtaking in mixed traffic using new variables. Finally, Section VI presents a discussion and conclusions.
\section{Methodology}
\subsection{Derivation of continuum model using the principle of area conservation}

\citet{lighthill1955kinematic}, and \citet{richards1956shock} developed the LWR model that describes traffic flow as a continuum, with traffic density and flow rate varying continuously in space and time. In mixed traffic, one can visualize vehicular traffic as a stream of continuum flow of different sizes of particles. This leads to deriving a continuum equation based on the principle of total vehicle area conservation, assuming:  
 \begin{enumerate}
     \item Traffic can be treated as a continuous medium using continuous density and flow-like variables.
     \item Total projected vehicles' area is conserved with suitable sinks and sources.
     \item A certain constitutive relationship holds between the areal density, areal flow, and speed. These relations may be theoretically based on a physical understanding but could be validated using empirical data. 
     \item Areal density and areal flow are set as a function of time-space.
     \item The variables are bound by the geometric and vehicle characteristics such as the road width, vehicle size, etc.
     \item Homogeneity need not hold at all points in space in any direction.
 \end{enumerate}

\subsubsection{Formulation}
\begin{figure}[h]
  \centering
  \includegraphics[width=1\textwidth]{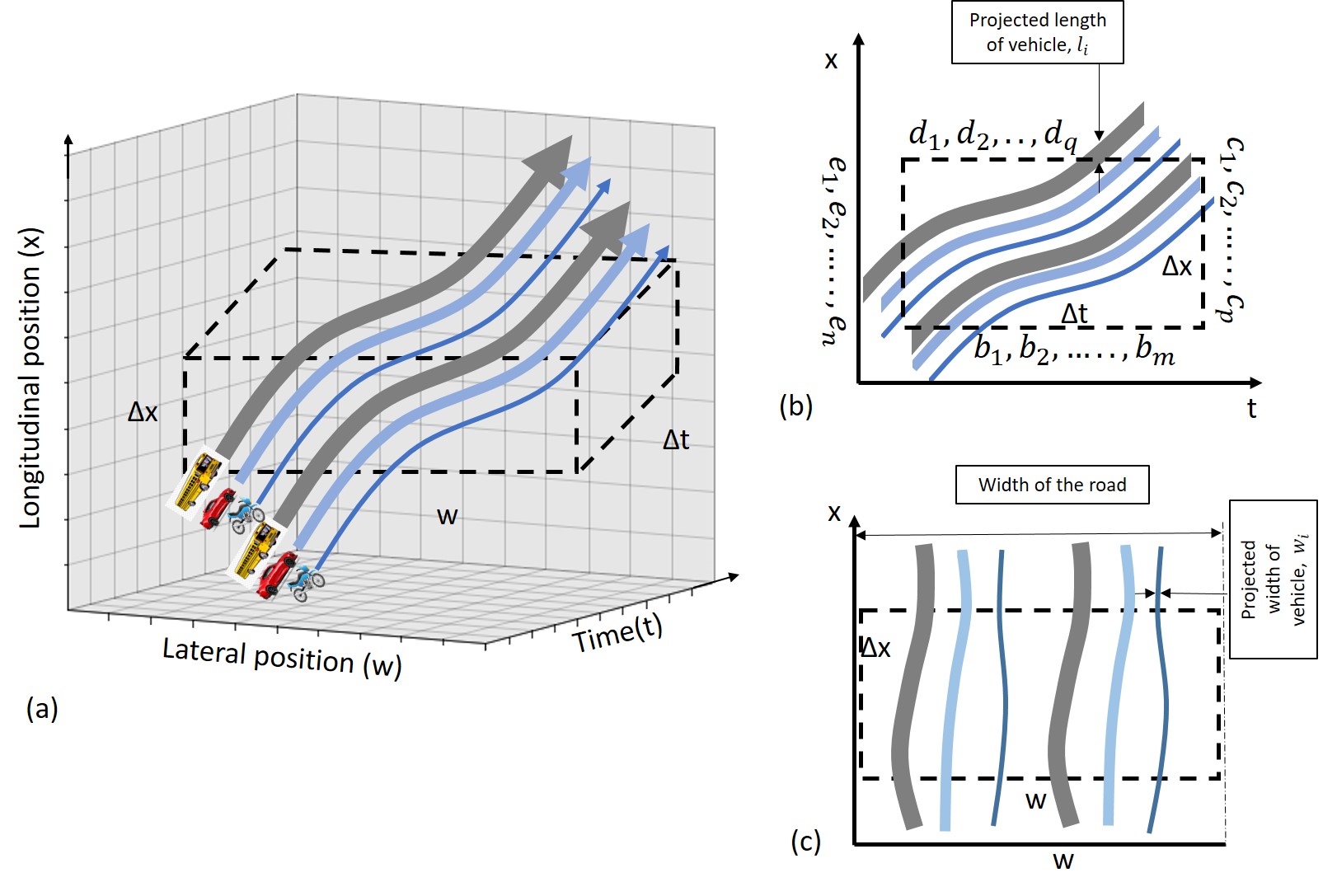}
  \caption{(a) Three-dimensional representation of time-space domain 
 and vehicle trajectories within the volume of $(\Delta x \times \Delta t \times w)$ (b) Two-dimensional representation of trajectories with size and number of vehicles coming in and going out in the domain of $(\Delta x \times \Delta t)$  (c) Two-dimensional representation of trajectories in the domain of $(\Delta x \times w)$ (color and width of the trajectories represents different physical characteristics of vehicles)}\label{fig:continuum}
\end{figure}

This section delves into formulating a continuum equation specifically addressing mixed traffic conditions. Figure \ref{fig:continuum} (a,b) shows the vehicle trajectories in a three-dimensional and two-dimensional representation within the considered volume of the region. Let $e_i,b_i,c_i,d_i$ be the area of the vehicle for different vehicle classes on various boundaries, as shown in Figure \ref{fig:continuum} (b). Therefore, applying the principle of vehicle area conservation, the total vehicle area entering the section will equal the total vehicle area exiting \eqref{eq:cnt_eqn0}, where source ($\delta_{source}$) and sink ($\delta_{sink}$) represent new vehicle area introduced or removed within the domain under consideration. Eq. \eqref{eq:cnt_eqn1} is derived by dividing $(\Delta x \Delta t w)$ to \eqref{eq:cnt_eqn0}.

\begin{subequations}\label{eq:litdiff}
\begin{align}
\sum_{i=1}^n e_i + \sum_{i=1}^m b_i+\delta_{source} &= \sum_{j=1}^p c_j + \sum_{j=1}^q d_j + \delta_{sink}\label{eq:cnt_eqn0}\\
\frac{\sum_{i=1}^n e_i-\sum_{j=1}^p c_j}{\Delta x \Delta t w} & = - \frac{\sum_{i=1}^m b_i-\sum_{j=1}^q d_j}{\Delta x \Delta t w} + \frac{\delta_{sink}-\delta_{source}}{\Delta x \Delta t w} \label{eq:cnt_eqn1}\\
\frac{k_a^e-k_a^c}{ \Delta t} & = - \frac{q_a^b-q_a^d}{\Delta x} + \frac{\delta_{sink}-\delta_{source}}{\Delta x \Delta t w} \label{eq:cnt_eqn2}\\
\frac{\Delta k_a}{ \Delta t} &= -\frac{\Delta q_a}{\Delta x }+f(x,t) \label{eq:cnt_eqn3}
\end{align}
\end{subequations}

 where $k_a^e = \frac{\sum_{i=1}^n e_i}{ \Delta x w} , k_a^c = \frac{\sum_{j=1}^p c_j}{ \Delta x w};  \hspace{2mm} q_a^b = \frac{\sum_{i=1}^m b_i}{ \Delta t w},  q_a^d = \frac{\sum_{j=1}^q d_j}{ \Delta t w}$.

In Equation \eqref{eq:cnt_eqn1}, the term $(\sum_{i=1}^n e_i - \sum_{j=1}^p c_j)$ represents the variation in vehicle area within the spatial segment $\Delta x$, while $(\sum_{i=1}^m b_i - \sum_{j=1}^q d_j)$ characterizes the change in vehicle area arising from vehicles traversing the boundaries over the interval $\Delta t$. The quantity $(\delta_{sink}-\delta_{source})$ denotes the net change in the vehicle area within the region $\Delta x$ $\Delta t$, representing the influence of sources and sinks within the region.

In Eq. \eqref{eq:cnt_eqn3}, $\Delta q_a$ is the rate of change in vehicle area over time and width of the roadway, and $\Delta k_a$ is the rate of change in vehicle area over the length and width of the roadway. From \eqref{eq:cnt_eqn3}, $f(x,t)$ represents the net flow of area per unit length and width of road per unit time. Assuming continuous medium of traffic flow, applying $lim_{t,x \rightarrow 0}$ of the Eq \eqref{eq:cnt_eqn3}.
\begin{equation}\label{eq:final_cont}
     \frac{\partial q_{a}}{\partial x } +\frac{\partial k_{a}}{ \partial t}=f(x,t)
\end{equation}
Similar to that in the LWR model, the areal continuum equation \eqref{eq:final_cont} takes a macroscopic approach by conserving vehicle area to model mixed traffic, treating traffic as a continuum fluid-like medium. The hyperbolic character of Equation \eqref{eq:final_cont} becomes apparent due to its incorporation of both spatial and temporal derivatives and describes traffic wave propagation like the LWR model. To effectively apply the areal continuum equation, it is essential to provide precise and well-defined definitions for the newly introduced areal variables, that are independent of the measurement method and traffic state.   
\subsubsection{Generalized Definition of Areal Traffic Flow Variables}
The fundamental variables, derived from the vehicle area conservation continuum model defined in Eq.\eqref{eq:final_cont}, are termed as areal flow $(q_a)$, areal density $(k_a)$, and areal speed $(v_a)$. Here, $q_a$ is defined as the total vehicle area crossing the road width at a location during the observation period. Similarly, $k_a$ is defined as the total vehicle area divided by the space length and road width. Areal speed represents the rate at which the road space is released.
\newline\newline \underline{Areal flow $(q_a)$}: Considering all the vehicles crossing at the entry of the road section (see Figure \ref{fig:continuum} (b)) during $\Delta t$ time period, where the vehicle area is represented by $b_i$ and $i$ is the index of vehicle, which total $m$ vehicles. We can define the variable $q_a$ according to the  \eqref{eq:cnt_eqn2}. Now by multiplying $(\Delta x \times m)$ in both denominator and numerator, we get the average area of the vehicles at the entry $(\frac{\sum_{i=1}^m b_i}{m}=\Bar{b})$, and total distance traveled by all the vehicles $(m \times \Delta x=\Bar{b} \sum dx_i)$. Therefore, $q_a$ can be defined as follows in \eqref{eq:qa} with the unit of $(m^2/hr-m)$.
\begin{equation}
\begin{split}\label{eq:qa}
     q_a & =\frac{\sum_{i=1}^m b_i}{ \Delta t \times w}= \frac{\sum_{i=1}^m b_i \times \Delta x \times m}{\Delta t \times w \times \Delta x \times m} = \frac{\sum_{i=1}^m b_i}{m} \times \frac{m \times \Delta x}{\Delta t \times \Delta x \times w}\\
     &=\frac{\Bar{b} \sum dx_i}{\Delta t \times \Delta x \times w}\\
     &=\frac{Average \hspace{1mm} vehicle \hspace{1mm} area  \times total \hspace{1mm} distance  \hspace{1mm} traveled}{Volume \hspace{1mm} of \hspace{1mm} the \hspace{1mm} region}.\\
    unit & \rightarrow \frac{m^2}{hr-m}
\end{split}
\end{equation}
\newline \underline{Areal density $(k_a)$}: Similarly to $q_a$, considering all vehicles in the $\Delta x$ road section at the initial time instant (see Figure \ref{fig:continuum} (b)), $k_a$ can be defined as follows in \eqref{eq:ka} with the unit of $(m^2/km-m)$. 
\begin{equation}
\begin{split}\label{eq:ka}
     k_a & =\frac{\sum_{i=1}^n e_i}{ \Delta x \times w}= \frac{\sum_{i=1}^n e_i \times \Delta t \times n}{\Delta x \times w \times \Delta t  \times n}= \frac{\sum_{i=1}^n e_i}{n} \times \frac{n\times \Delta t}{\Delta t \times \Delta x \times w}\\
     &=\frac{\Bar{e} \sum dt_i}{\Delta t \times \Delta x \times w}\\
     &=\frac{Average \hspace{1mm} vehicle \hspace{1mm} area  \times total \hspace{1mm} time  \hspace{1mm} spent}{Volume \hspace{1mm} of \hspace{1mm} the \hspace{1mm} region}.\\
     unit & \rightarrow \frac{m^2}{km-m}
\end{split}
\end{equation}
\newline \underline{Areal speed $(v_a)$}: Rate of vehicle area released in $(\Delta x \times \Delta t)$ region, which is given as the ratio of the $q_a$ and $k_a$. Since the average vehicle area during the observation period remains the same in the region, 
\begin{equation}\label{eq:va}
    \begin{split}
        v_a &=\frac{q_a}{k_a}
        =\frac{ \sum dx_i}{ \sum dt_i}\\
        &=\frac{total \hspace{1mm} distance  \hspace{1mm} traveled}{total \hspace{1mm} time  \hspace{1mm} spent}.\\
     unit & \rightarrow \frac{km}{hr}
    \end{split}
\end{equation}
Note that the definition of areal speed is the same as the space mean speed from Edie's generalized definitions \citep{edie1963discussion}. Above are the generalized definitions that are independent of the measurement method. These explanations appear unfamiliar when considering our broad understanding of flow and density. Still, they align with prior definitions, regardless of specific measurement techniques in their approach to averaging. 
\par
\begin{figure}[h]
\centering  
\subfigure[]{\includegraphics[width=0.48\linewidth]{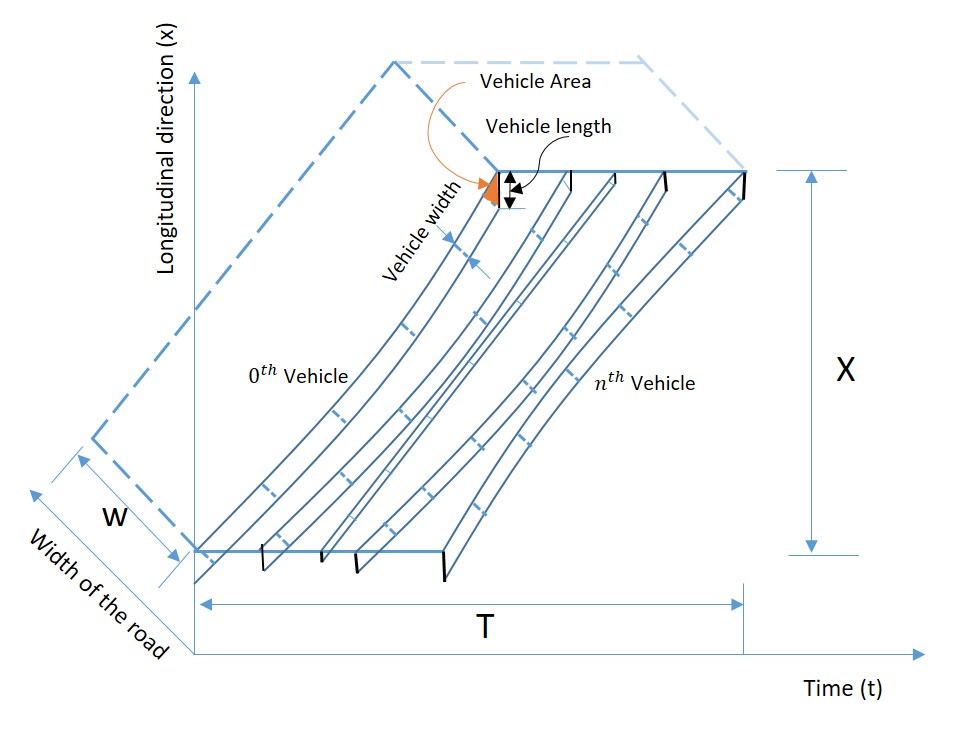}}
\subfigure[]{\includegraphics[width=0.48\linewidth]{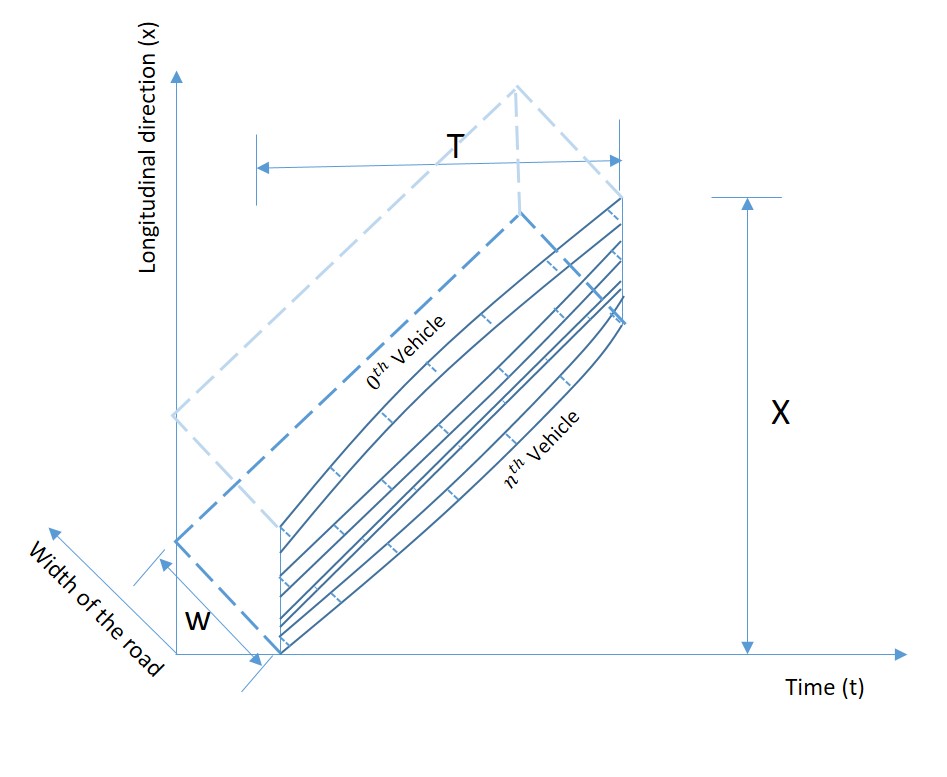}}
\caption{ Trajectories of a group of n vehicles over (a) a distance X, and (b) over a time T }
\label{fig:edie_def}
\end{figure}

 A theoretical derivation shows that the proposed definitions for areal flow and density can be obtained from the generalized definitions. Figure \ref{fig:edie_def} (a) shows a time-space domain of vehicles $0^{th}$ to $N^{th}$ trajectories in a section of a roadway of length $X$. The analysis is presented within the spatio-temporal domain. Here, the projected width and length of the road are shown, along with the vehicle area. Dotted and continuous lines around the vehicle trajectories depict the vehicle dimensions (width and length), which are aligned with the road's width and longitudinal direction, respectively. By assuming the linear slope of the trajectory of $0^{th}$ and $N^{th}$ vehicles as $v_0$ and $v_N$, the volume of the considered region is expressed as \eqref{eq:ed_vol}.

\begin{equation}\label{eq:ed_vol}
    Vol.=TXW-\frac{X^2W}{2}\left( \frac{1}{v_0}+\frac{1}{v_N}\right)
\end{equation}
Where $T$ is the total time of observation and $v_0,v_N$ is the average speed in the longitudinal direction for the two boundary vehicles.
\begin{equation}\label{eq:ed_q}
    \begin{split}
        q_a &=\frac{Average \hspace{2mm} vehicle \hspace{2mm} area \times total \hspace{2mm} distance  \hspace{2mm} travel}{Volume \hspace{2mm} of \hspace{2mm} the \hspace{2mm} region}\\
        &= lim_{X \rightarrow \Delta x }  \frac{\left(\frac{1}{N}\sum_{i=0}^N a_i\right) NX}{TXW-\frac{X^2W}{2}\left( \frac{1}{v_0}+\frac{1}{v_N}\right)} \approx \frac{\left(\frac{1}{N}\sum_{i=0}^N a_i\right) \sum_{i=0}^N \Delta x}{T\Delta x W}
    \end{split}
\end{equation}
Here $a_i$ represents the projected area of the $i^{th}$ vehicle. Eq. \eqref{eq:ed_q} reduces to \eqref{eq:qa}, $q_a= \frac{ \Bar{a} \sum \Delta x}{T\Delta x W}$, by considering measurement at a very small distance by $ \Delta x $. The denominator term of \eqref{eq:ed_q}, $ -\frac{\Delta x^2W}{2}\left( \frac{1}{v_0}+\frac{1}{v_N}\right) \rightarrow 0$, become very small; hence this term can be ignored. Thus, \eqref{eq:ed_q} reduces to $q_a$ for measurement assumed to have been made over a very short distance. Since the paper assumes a finite area of the vehicle, $\Delta x$ cannot tend to zero. Similarly, $k_a$ can be expressed as \eqref{eq:ed_kka}. 
\begin{equation}\label{eq:ed_kka}
    \begin{split}
        k_a &=\frac{Average \hspace{2mm} vehicle \hspace{2mm} area \times total \hspace{2mm} time  \hspace{2mm} spent}{Volume \hspace{2mm} of \hspace{2mm} the \hspace{2mm} region}\\
        &= lim_{X \rightarrow \Delta x } \frac{\left(\frac{1}{N}\sum_{i=0}^N a_i\right) \sum_{i=0}^N \frac{X}{v_i}}{TXW-\frac{X^2W}{2}\left( \frac{1}{v_0}+\frac{1}{v_N}\right)} \approx \frac{\left(\frac{1}{N}\sum_{i=0}^N a_i\right) \sum_{i=0}^N \frac{\Delta x}{v_i}}{T\Delta x W}
    \end{split}
\end{equation}
Which becomes \eqref{eq:ka}, $k_a= \frac{\left(\frac{1}{N}\sum_{i=0}^N a_i\right)\sum_{i=0}^N \frac{\Delta x}{v_i}}{T\Delta x W}$ as $X$ approaches $\Delta x$ and the size of the measurement length becomes very small. Thus, unlike Edie's definition, the flow exhibited by such a platoon remains affected by the area of intermediate vehicles and can not be solely determined by the first and last vehicles.
\par Figure \ref{fig:edie_def} (b) shows a time-space region enclosed by a period T with the vehicle trajectories $0$ to $N$. The volume and stream parameters described in Figure \ref{fig:edie_def} (b) are as follows.

\begin{subequations}\label{eq:ed_vol1}
\begin{align}
    Vol. &=TXW-\frac{T^2W}{2}\left( \frac{1}{v_0}+\frac{1}{v_N}\right)\\
    q_a &= lim_{T \rightarrow \Delta t } \frac{\left(\frac{1}{N}\sum_{i=0}^N a_i\right) \sum_{i=0}^N T.v_i}{TXW-\frac{T^2W}{2}\left( \frac{1}{v_0}+\frac{1}{v_N}\right)} \approx \frac{\left(\frac{1}{N}\sum_{i=0}^N a_i\right) \sum_{i=0}^N \Delta t v_i}{X\Delta tW} \\
    k_a &= lim_{T \rightarrow \Delta t } \frac{\left(\frac{1}{N}\sum_{i=0}^N a_i\right) NT}{TXW-\frac{T^2W}{2}\left( \frac{1}{v_0}+\frac{1}{v_N}\right)} \approx \frac{\left(\frac{1}{N}\sum_{i=0}^N a_i\right) \sum \Delta t}{XW\Delta t}
\end{align}
\end{subequations}
Note that as $T$ approaches a very small value $\Delta t$, \eqref{eq:ed_vol1} (b) and \eqref{eq:ed_vol1} (c) reduces to the definitions of areal flow in \eqref{eq:qa} and areal density in \eqref{eq:ka}. These mathematical expressions can be applied in scenarios where T is substantially larger than zero, facilitating the computation of mean flow magnitudes present within a platoon of $N$ vehicles over the provided duration.

Therefore, the generalized definitions of areal variables demonstrably hold true across all measurement frameworks. Also, note that these areal variables are particularly suitable for lane-free movement scenarios, as they take into account the entire width of the road rather than being confined to a single lane.

\subsubsection{Relationship among traffic flow variables}
This section investigates the relationships between the newly introduced areal variables and established traffic state variables like density $(k)$, flow $(q)$, occupancy $(O_c)$, and area occupancy $(ao)$. Investigating these interrelationships helps visualize the new variables in the traditional context but also provides valuable insights for transferability and interoperability across the areal variables and traditional variables. Furthermore, these relationships could be instrumental in characterizing traffic conditions within various analytical frameworks and potentially seamlessly utilizing existing theories for traffic analysis. The relationships between these variables are derived and shown in \eqref{eq:k_ka}-\eqref{eq:ka_oa1}. 
\newline
\newline \underline{Areal and Conventional variables:}
Using \eqref{eq:qa}-\eqref{eq:va} and Edie's generalized definitions \citep{edie1963discussion},
\begin{subequations}\label{eq:k_ka}
\begin{align}
k_a&=\frac{\sum_{i=1}^n e_i}{\Delta x w}
        =\frac{ \Bar{e} \sum dt_i}{\Delta t \times \Delta x \times w}
        =\frac{ \Bar{e}} {w}k
         \label{eq:k_ka_1}\\
q_a & =\frac{\sum_{i=1}^m b_i}{ \Delta t w}
         =\frac{ \Bar{b} \sum dx_i}{\Delta t \times \Delta x \times w}
        =\frac{ \Bar{b}} {w}q
         \label{eq:q_qa_1}\\
v_a &=\frac{q_a}{k_a}=\frac{q}{k}=v_s
\end{align}
\end{subequations}

Where $\Bar{e}$ and $\Bar{b}$ represent the average vehicle area. Thus, the areal density and areal flow are the scaled versions of the density and flow, respectively, where the scaling factor is the average vehicle area divided by the road width. Based on the average vehicle dimensions and road width, it can be observed that generally, the areal density and areal flow are higher than the conventional density and flow values. Additionally, for two similar roads with similar $k$ and $q$, the \eqref{eq:k_ka_1}, and \eqref{eq:q_qa_1} would mean that the road with a higher proportion of bigger vehicles has higher $k_a$ and $q_a$. Also, across two roads with varying widths, the wider road will exhibit lower areal density and areal flow values compared to the narrower road.
\newline
\newline \underline{Occupancy $(O_c)$ and Areal Density $(k_a)$:}

For mixed traffic with class-specific length $(L_i)$ and width $(B_i)$, crossing over a detector with length $(d)$, the $O_c$ can be defined as follows:
\begin{equation}\label{eq:ka_oc1}
    \begin{split}
        O_c&=\frac{\sum_i (L_i+d)/v_i}{T}\\
        &=\frac{N}{T}\frac{1}{N}\sum \frac{L_i}{v_i}+d\frac{N}{T}\frac{1}{N}\sum \frac{1}{v_i}
        =\frac{q}{N}\sum \frac{L_i}{v_i}+d.q.\frac{1}{N}\sum \frac{1}{v_i}\\
        &=\frac{q}{v_s} \frac{v_s}{N}\sum \frac{L_i}{v_i}+d.\frac{q}{v_s}\\
        &=k \frac{1}{1/N \sum 1/v_i}\frac{1}{N}\sum \frac{L_i}{v_i}+d.k\\
        &=k \frac{\sum \frac{L_i}{v_i}}{\sum \frac{1}{v_i}}+d.k 
        =k_a\frac{w}{\frac{1}{N} \sum L_i B_i}( \frac{\sum \frac{L_i}{v_i}}{\sum \frac{1}{v_i}}+d) 
        =\frac{k_awN}{\sum L_i B_i}( \frac{\sum \frac{L_i}{v_i}}{\sum \frac{1}{v_i}}+d)\\
    \end{split}
\end{equation}
For homogeneous conditions $L_i=L$ and $B_i=B$. Therefore,
\begin{equation}\label{eq:ka_oc2}
    \begin{split}
        O_c&=\frac{k_aw}{L B}( L+d)\\
        &=k_a \frac{w}{B}(1+\frac{d}{L})\\
    \end{split}
\end{equation}
Thus, in mixed traffic, $O_c$ is lower on roads with a higher proportion of trucks but higher on wider roads.
\newline 
\newline \underline{Area Occupancy $(ao)$ and Areal Density $(k_a)$:}
For mixed traffic, $ao$ can be defined as follows: 
\begin{equation}\label{eq:ka_oa1}
    \begin{split}
       ao&=\frac{\sum \frac{ (L_i+d)B_i}{v_i}}{Tw}\\
       &= \frac{\sum \frac{ (L_i+d)B_i}{v_i} \frac{1}{N}}{\frac{Tw}{N}}
       =\frac{\sum \frac{ (L_i+d)B_i}{v_i} \frac{q}{N}}{w} \\
       &=\frac{\sum \frac{ (L_i+d)B_i}{v_i} k\frac{v_s}{N}}{w} 
       =\frac{\sum \frac{ (L_i+d)B_i}{v_i} \frac{k}{\sum \frac{1}{v_i}}}{w}\\
       &=\frac{k}{w}\frac{\sum \frac{ (L_i+d)B_i}{v_i}}{\sum \frac{1}{v_i}}
       =\frac{k_a}{\frac{1}{N} \sum L_i B_i}\frac{\sum \frac{ (L_i+d)B_i}{v_i}}{\sum \frac{1}{v_i}}\\
       &=\frac{k_aN}{\sum L_i B_i}(\frac{\sum \frac{L_iB_i}{v_i}}{\sum \frac{1}{v_i}}+\frac{d\sum \frac{B_i}{v_i}}{\sum \frac{1}{v_i}})
    \end{split}
\end{equation}

The area occupancy in mixed traffic is solely determined by the weighted average area of vehicles and the individual sizes of those vehicles. In contrast to regular occupancy, it is not influenced by the road's geometry. For homogeneous traffic, \eqref{eq:ka_oa1} reduces to $ao=k_a(1+\frac{d}{L})$. Thus, the scaling factor for \textit{ao} of homogeneous traffic is a constant value, a power function of $L$, solely dependent on the length of the vehicle without any other variables. This indicates that $ao$ is lower on roads with a higher proportion of large vehicles. Also, the multiplying factor of over 1 indicates that $ao$ is always greater than $k_a$ (note that both these are dimensionless). Note that in homogeneous traffic, $O_c$ vs. $k$ is analogous to $ao$ vs. $k_a$ in the sense that the multiplying factor is only a function of $L$ and $d$ and is independent of the vehicle width or road width.
\par
\begin{figure}[h]
\centering  
\subfigure[$k_a-k$]{\includegraphics[width=0.32\linewidth]{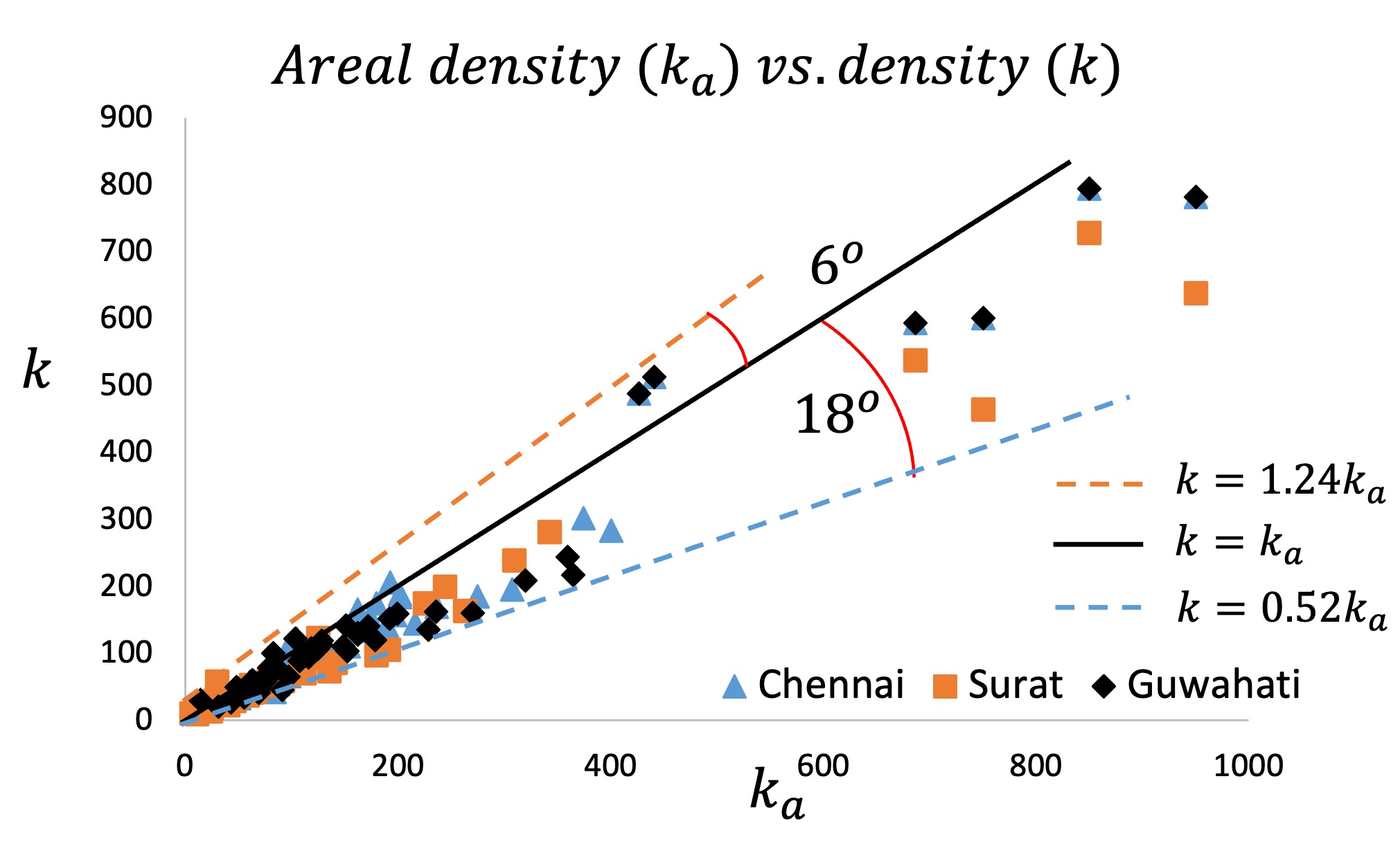}}
\subfigure[$k_a-O_c$]{\includegraphics[width=0.32\linewidth]{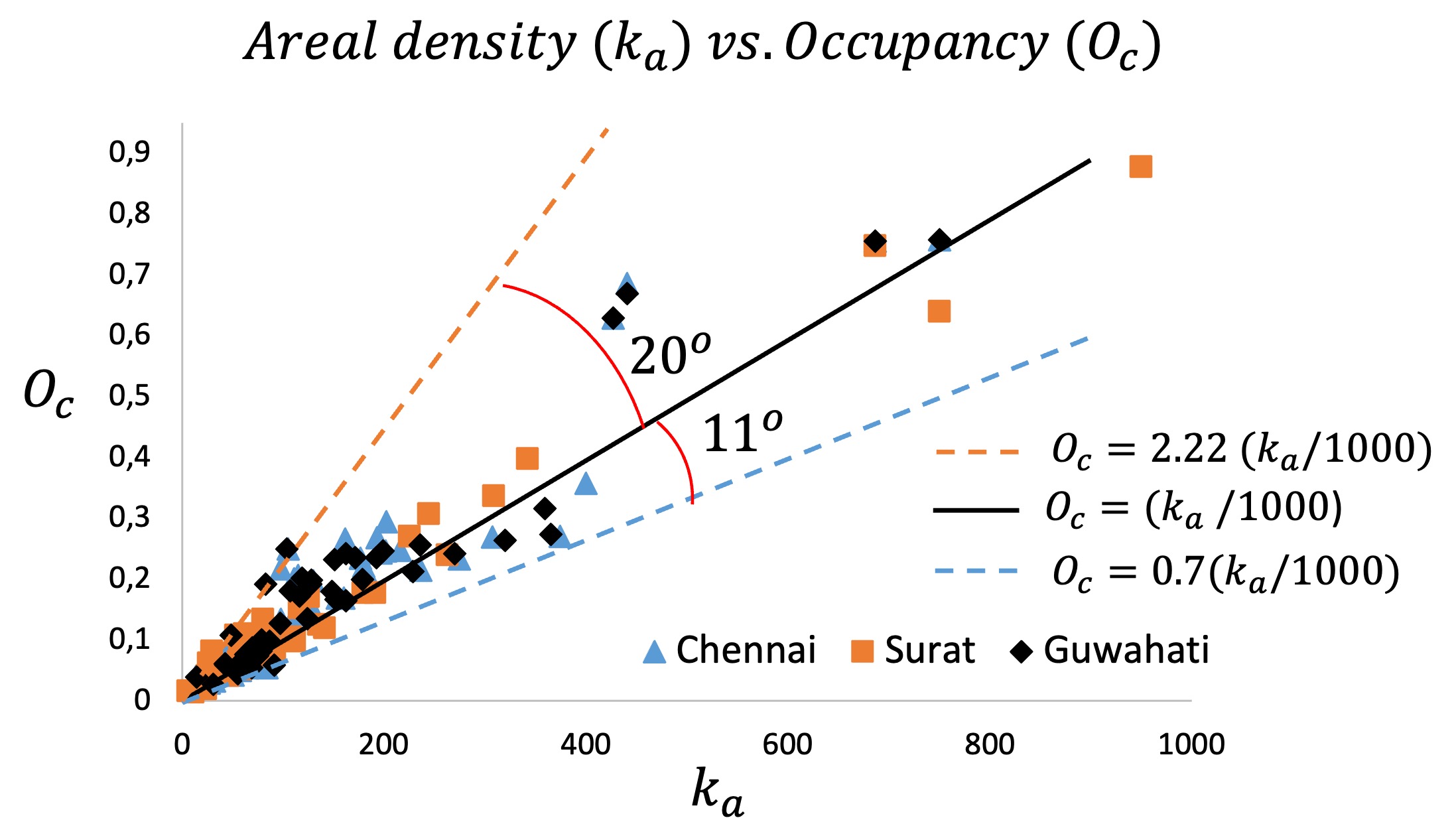}}
\subfigure[$k_a-ao$]{\includegraphics[width=0.32\linewidth]{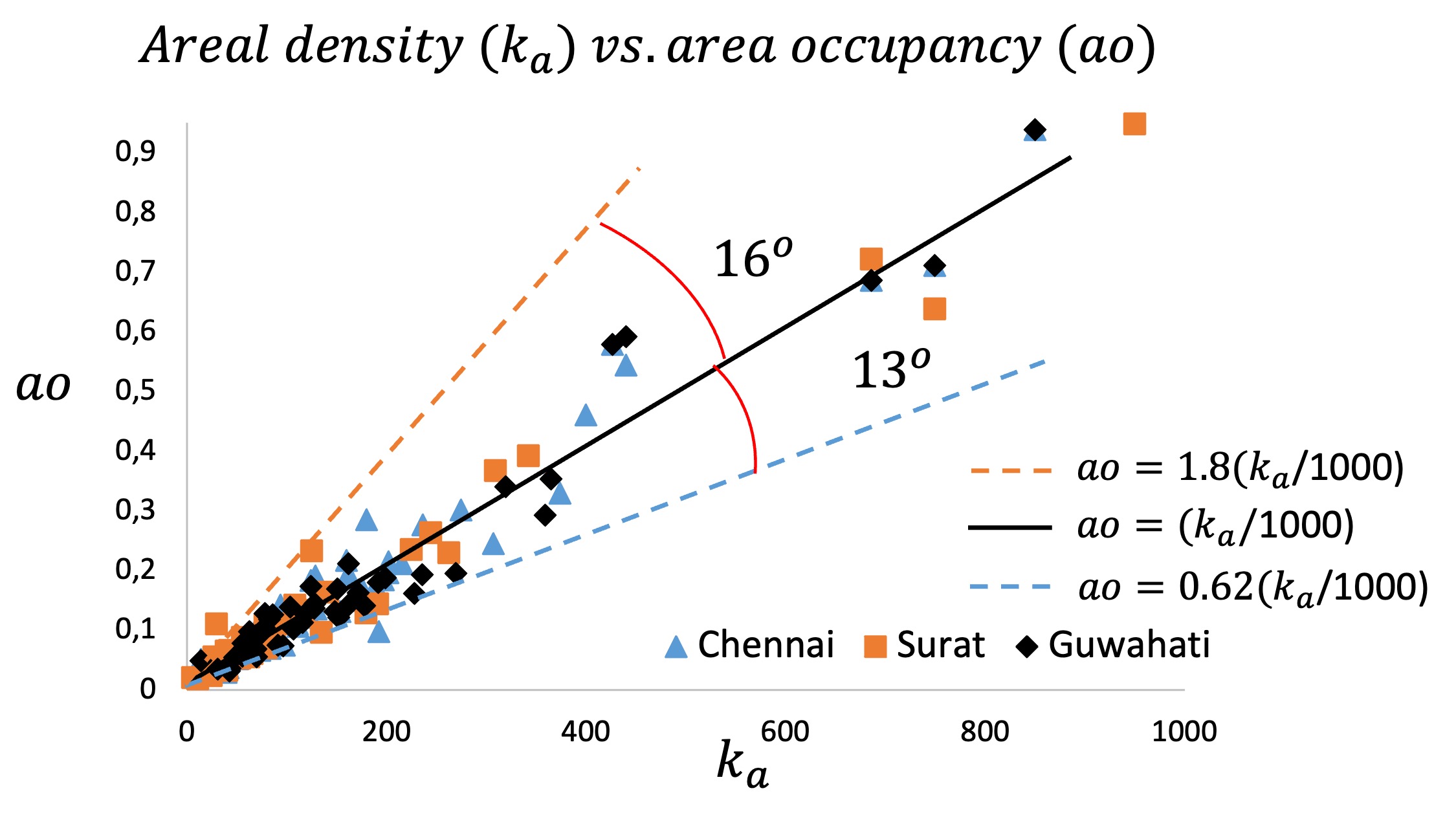}}
\caption{Empirical relationships between the proposed areal density $(k_a)$ with (a) density $(k)$, (b) occupancy $(O_c)$, and (c) area occupancy $(ao)$.}
\label{fig:relationships}
\end{figure}

 Further analysis was conducted to explore the relationships between areal density ($k_a$) and other critical traffic parameters using empirical datasets. Figure \ref{fig:relationships} delineates a linear increase in $k$, $O_c$, and $ao$ as functions of $k_a$, suggesting a systematic scaling of these parameters with $k_a$. The inset equations in Figure \ref{fig:relationships} establish the upper-bound (dotted orange line) and lower-bound (dotted blue line) linear relationships that frame the dispersion of data points across the measured variables. The one-to-one correlation is effectively illustrated by a continuous black line, which acts as a comparative baseline against the observed data trends. Figure \ref{fig:relationships} (a) demonstrates that $k_a$ predominantly resides below the density ($k$), capturing the physical size variances of vehicles in mixed traffic compositions. This observation supports the premise that density, typically measured by vehicle counts, may not fully account for the physical space occupied by larger vehicles. Given the minimal presence of HVs in the traffic mixes studied, their larger size disproportionately influences the areal density, justifying the placement of data points below the one-to-one baseline in accordance with Eq. \eqref{eq:k_ka}a. In Figure \ref{fig:relationships} (b), occupancy ($O_c$) appears slightly elevated, reflecting a reduced proportion of HVs, thus aligning with Eq. \eqref{eq:ka_oc1}. This indicates that $O_c$ increases at a rate slightly higher than $k_a$. Figure \ref{fig:relationships} (c) presents area occupancy ($ao$) slightly exceeding $k_a$, validating Eq. \eqref{eq:ka_oa1}. This analysis offers a more profound understanding of the real-world data supporting the theoretical explanation provided above.

\section{steady-state relationships among areal variables}
\subsection{Data description}
\begin{figure}[H]
\centering 
\centering
  \includegraphics[width=0.7\textwidth]{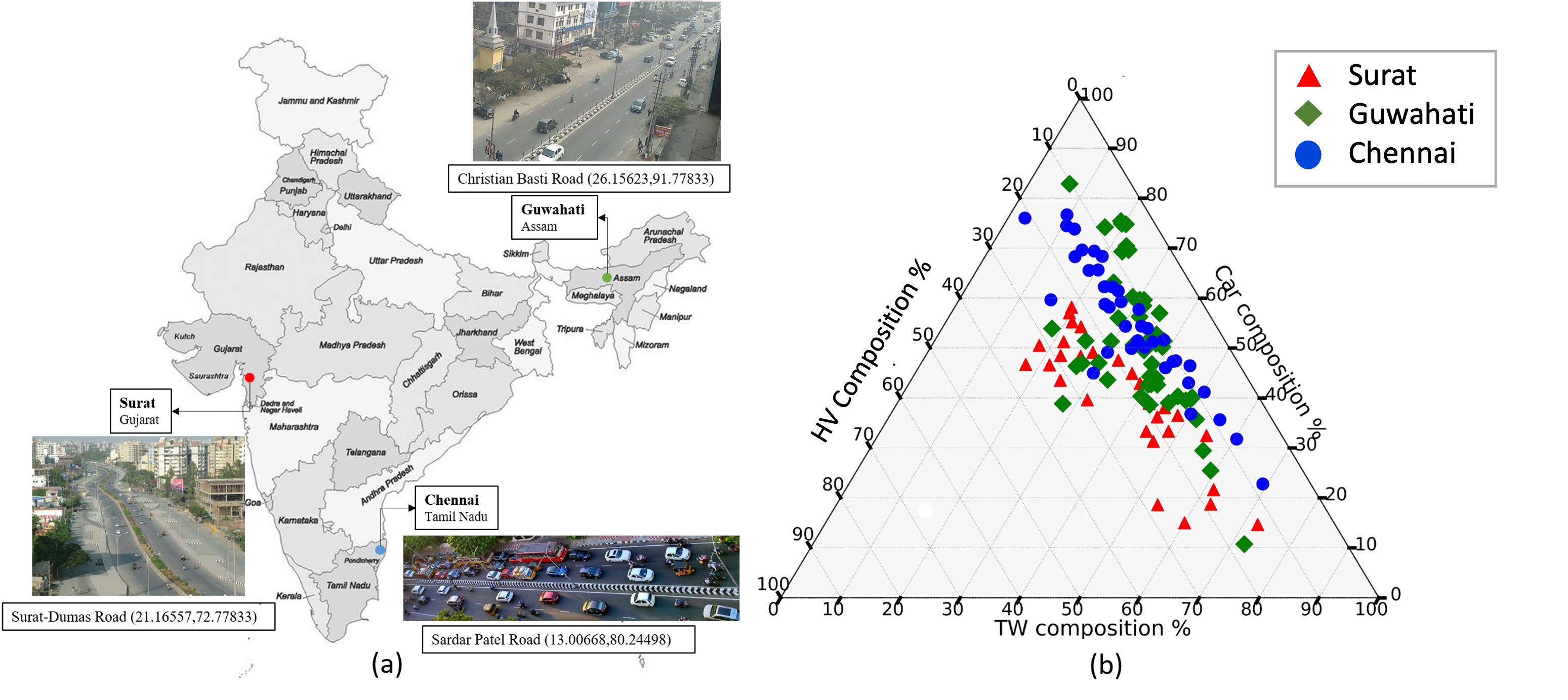}
  \caption{Data collection locations (a) Chennai, Lat-lon position of the road (13.00668,80.24498); Surat, Lat-lon position of the road (21.16557,91.77833); Guwahati, Lat-lon position of the road (26.15623,91.77833) (b) Ternary plot of percentage proportion of different vehicle classes in the stream at different locations \citep{maiti2024universality}. }
\label{fig:data}
\end{figure}

This section presents the steps undertaken to develop the fundamental diagrams using field data collected across the three cities in India, Chennai, Surat, and Guwahati.  Traffic data used in the analysis were collected from a six-lane divided carriageway  (width= 10.5 m in each direction ) of an intercity corridor under the mixed traffic environment. All three roads are major state highways, and the length of the segments varies from 100 m to 200 m. Data collection locations are the mid-blocks free from the influence of traffic signals, roundabouts, or cross streets. Figure \ref{fig:data} (a) shows the top view of the three data collection locations used in this paper. The vehicle trajectories were extracted from the video using the image processing techniques discussed in \cite{dhatbale2021deep,maiti2023empirical}. Observed vehicles are categorized into three groups based on physical dimensions: cars, two-wheelers (TWs), and heavy vehicles (buses, trucks). The percentage of cars varies from 45\% to 55\% for all locations, while TW share is 33\%-35\% of total traffic. The observed physical dimensions of the vehicles are presented in Table \ref{tab:dimen}.  The typical length of a vehicle varies between 1.8 m and 10.5 m between TWs and buses, while the widths of vehicles range from 0.6 m to 2.5 m. TWs and cars travel with a maximum speed of 65 km/hr, followed by ThW, LCV, buses, and trucks. Figure \ref{fig:data} (b) illustrates the compositions of TW, cars, and buses-trucks-LCVs (HVs) in the study locations. 

\begin{table}[h]
\centering
\caption{Vehicle characteristics in the study locations}
\label{tab:dimen}
\begin{tabular}{@{}ccccc@{}}
\toprule
Class & Length & Width & Max. observed  & Min. observed  \\ 
 &  &  & speed (Km/hr) &  speed(Km/hr) \\ \midrule
TW    & 1.8     & 0.6  & 65 & 5 \\
Car   & 4.7      & 1.7 &  65  & 3   \\
Truck   & 8.4      & 2.5  &   45 & 3   \\ 
Bus   & 10.5     & 2.5 &   50 &   3  \\ 
\bottomrule
\end{tabular}%
\end{table}

\subsection{Stream fundamental diagrams}

\begin{figure}[h]
\centering  
\subfigure[cumulative area and oblique plot of cumulative area curve ]{\includegraphics[width=0.45\linewidth]{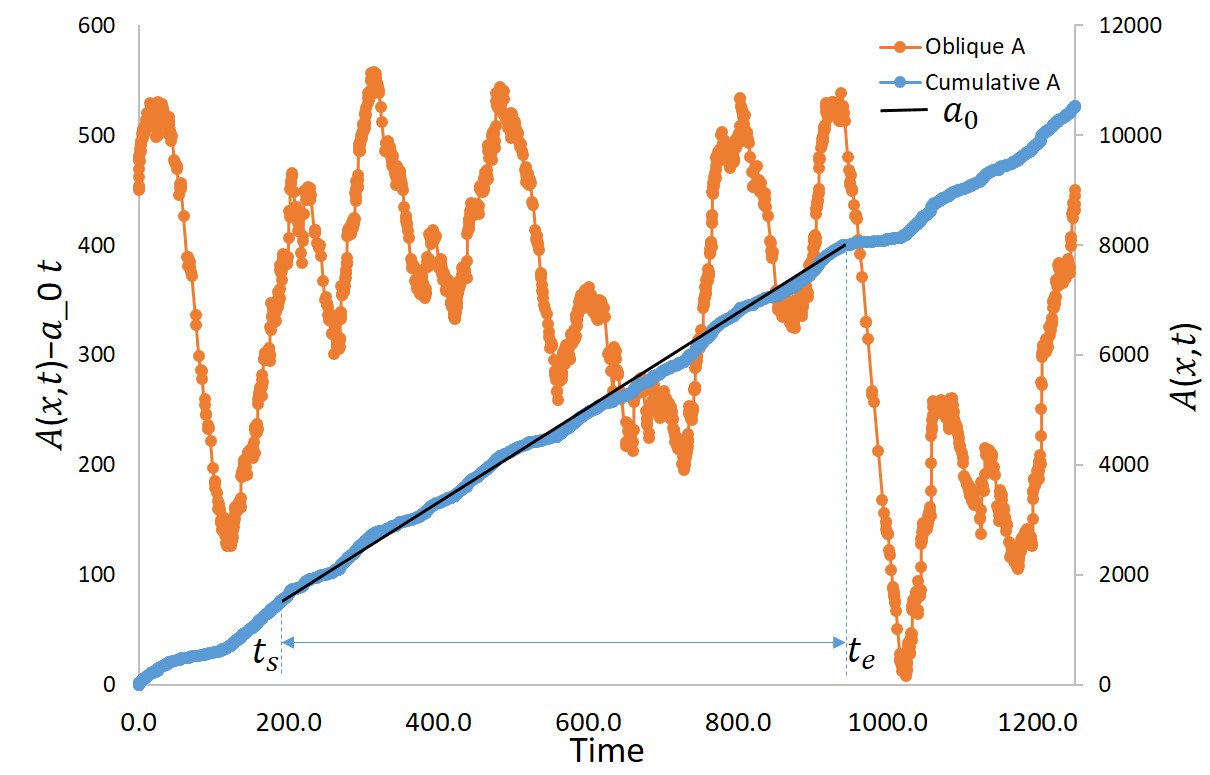}}
\subfigure[Oblique cumulative area curve and Oblique cumulative occupancy curve]{\includegraphics[width=0.45\linewidth]{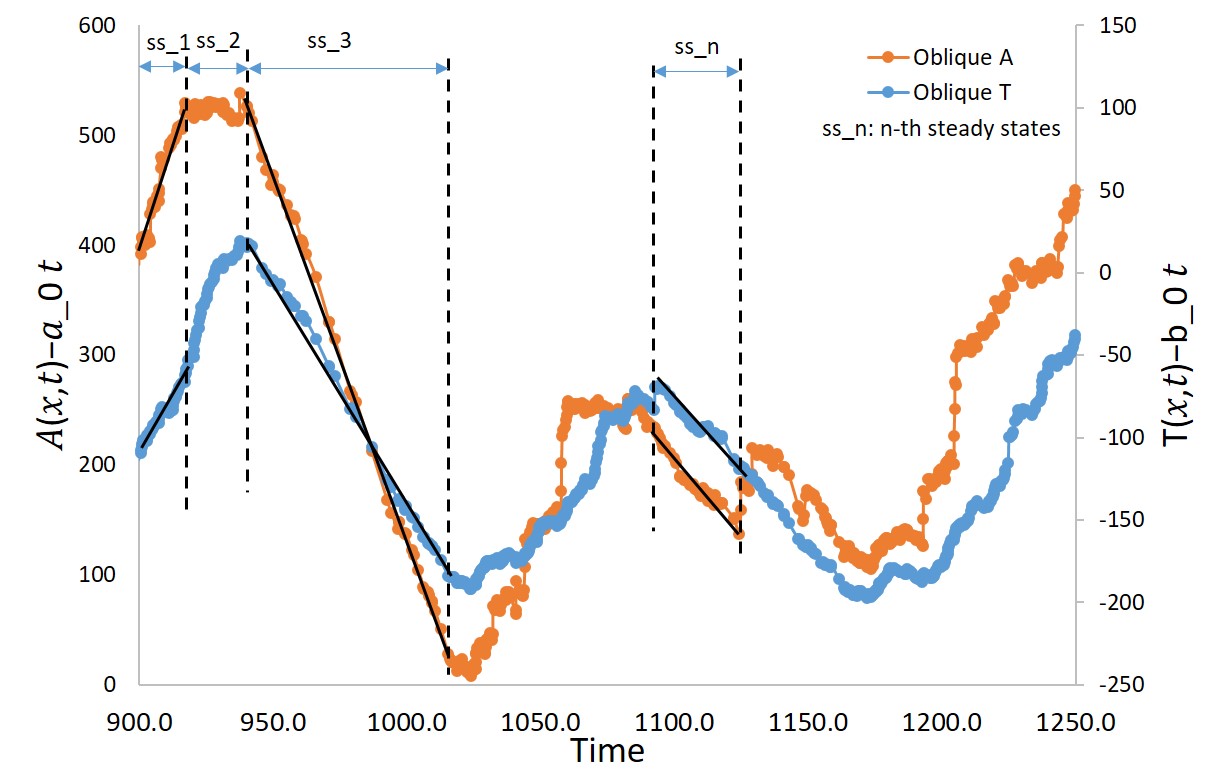}}
\caption{Steady-states identification using cumulative area and occupancy plots with the help of oblique cumulative vehicle area and occupancy plots.}
\label{fig:cac}
\end{figure}

This section explores the empirical bivariate relationships among the proposed areal variables. Using the approach outlined by \citet{cassidy1998bivariate}, we identified oblique plots and steady states. At a detector location, \( X \), we defined a cumulative vehicle area plot by summing the projected area of arriving vehicles over time, \( A(X, t) \), along with a cumulative occupancy plot, \( T(X, t) \). Figure \ref{fig:cac}(a) depicts a cumulative area curve and its corresponding oblique plot on the primary axis. The average area arrival rate was determined by identifying the start \((t_s)\) and end \((t_e)\) times of the linear trend in the cumulative area plot, calculated as \( a_0 = \frac{A(X, t_e) - A(X, t_s)}{t_e - t_s} \). Similarly, we estimated the average occupancy rate as \( b_0 = \frac{T(X, t_e) - T(X, t_s)}{t_e - t_s} \). Steady-state periods were defined as the intervals during which the oblique plots of cumulative vehicle area, \( A(X, t) - a_0t \), and cumulative occupancy time, \( T(X, t) - b_0t \), exhibit linear slopes. Figure \ref{fig:cac} (b) illustrates the common durations \((ss_1, ss_2, \ldots, ss_n)\) where both oblique plots are linear, indicating near-stationary states.

This study evaluates six established speed-density fundamental diagram models to determine their suitability for mixed traffic conditions using areal variables. The expressions of these models are presented in Table \ref{tab:speed_density_models}.
\begin{table}[H]
\centering
\caption{Summary of Speed-Density Relationships for Various Traffic Flow Models*}
\label{tab:speed_density_models}
\begin{tabular}{lc}
\hline
\textbf{Model} & \textbf{Speed-Density Function}  \\ \hline
\citet{greenshields1935study} & $v(k_a) = v_{max} \left( 1 - \dfrac{k_a}{k_{a,jam}} \right)$  \\
\citet{grinbeerg1959analysis}  & $v(k_a) = v_{crit} \ln \left( \dfrac{k_{a,jam}}{k_a} \right)$ \\
\citet{underwood1961speed} & $v(k_a) = v_{max} e^{-\frac{k_a}{k_{a,crit}}}$   \\
\citet{castillo1995functional}  & $v(k_a) = v_{max} \left[ 1 - \exp \left( 1 - \exp \left( \dfrac{\omega_a}{v_{max}} \left(\dfrac{k_{a,jam}}{k_a}-1\right) \right) \right) \right]$ \\
\citet{daganzo1994cell}  & 
        $v(k_a)= 
   \begin{cases}
    v_{max} & \text{if } 0\leq k_a \leq k_{a,crit}\\
    \omega_a\left(\frac{k_{a,jam}}{k_a}-1\right),              & \text{otherwise}
    \end{cases} $\\
\citet{smulders1990control} & 
$v(k_a) =   
\begin{cases}
v_{max}-\frac{(v_{max}-v_{crit})k_a}{k_{a,crit}} & \text{if } 0\leq k_a \leq k_{a,crit} \\
\omega_a \left(\dfrac{k_{a,jam}}{k_a}-1\right) & \text{otherwise}
\end{cases}$\\
\bottomrule
\end{tabular}
\end{table}

$*$ where stream wave speed $(\omega_a)$ is expressed as \eqref{eq:str_wave}.
\begin{equation}\label{eq:str_wave}
   \omega_a =\frac{v_{crit} k_{a,crit}}{k_{a,jam}-k_{a,crit}}
\end{equation}

\par
\begin{figure}[H]
\centering 
\subfigure[Chennai-Smulders]{\includegraphics[width=0.3\linewidth]{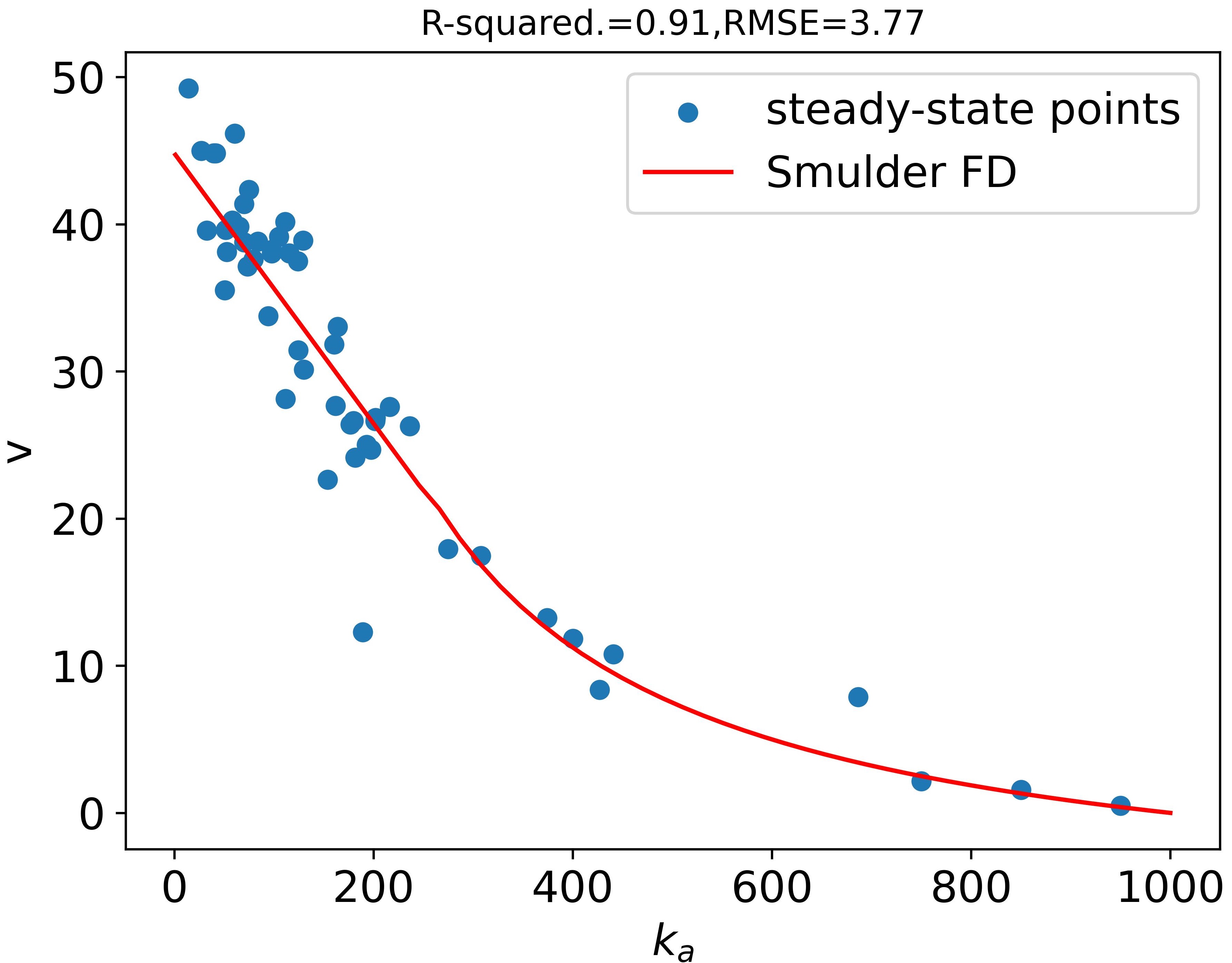}}
\subfigure[Surat-Smulders]{\includegraphics[width=0.3\linewidth]{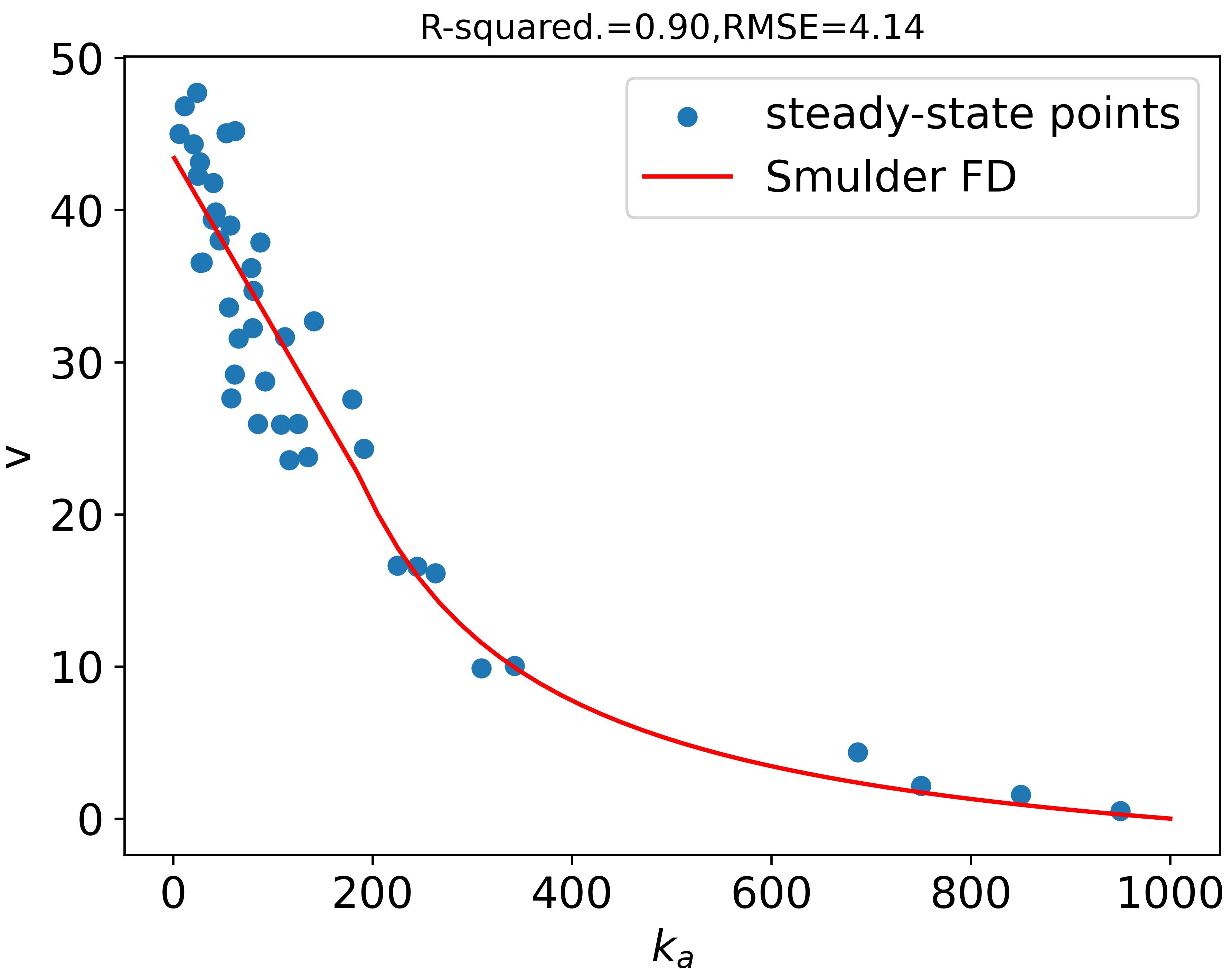}}
\subfigure[Guwahati-Smulders]{\includegraphics[width=0.3\linewidth]{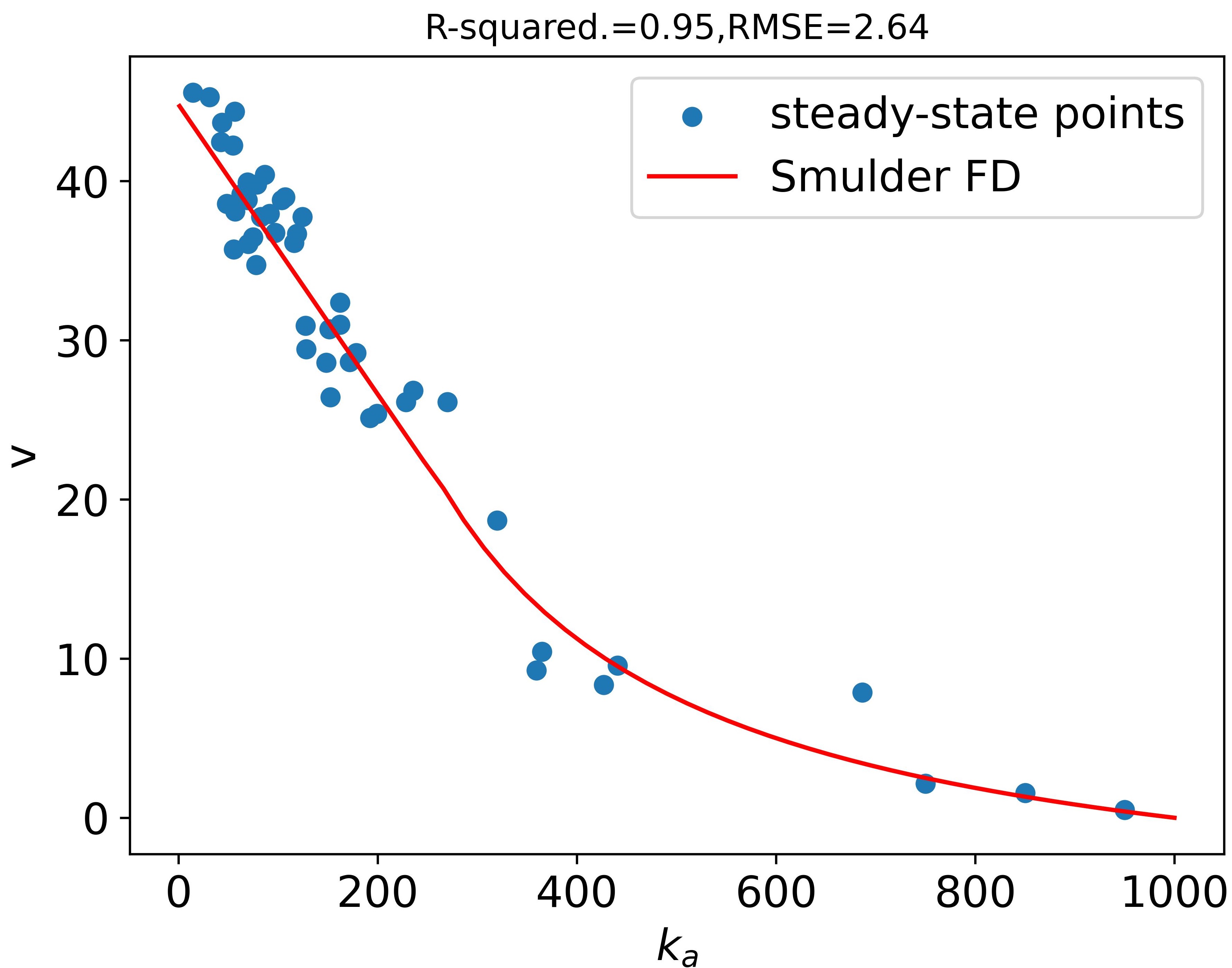}}
\subfigure[Chennai-Smulders]{\includegraphics[width=0.3\linewidth]{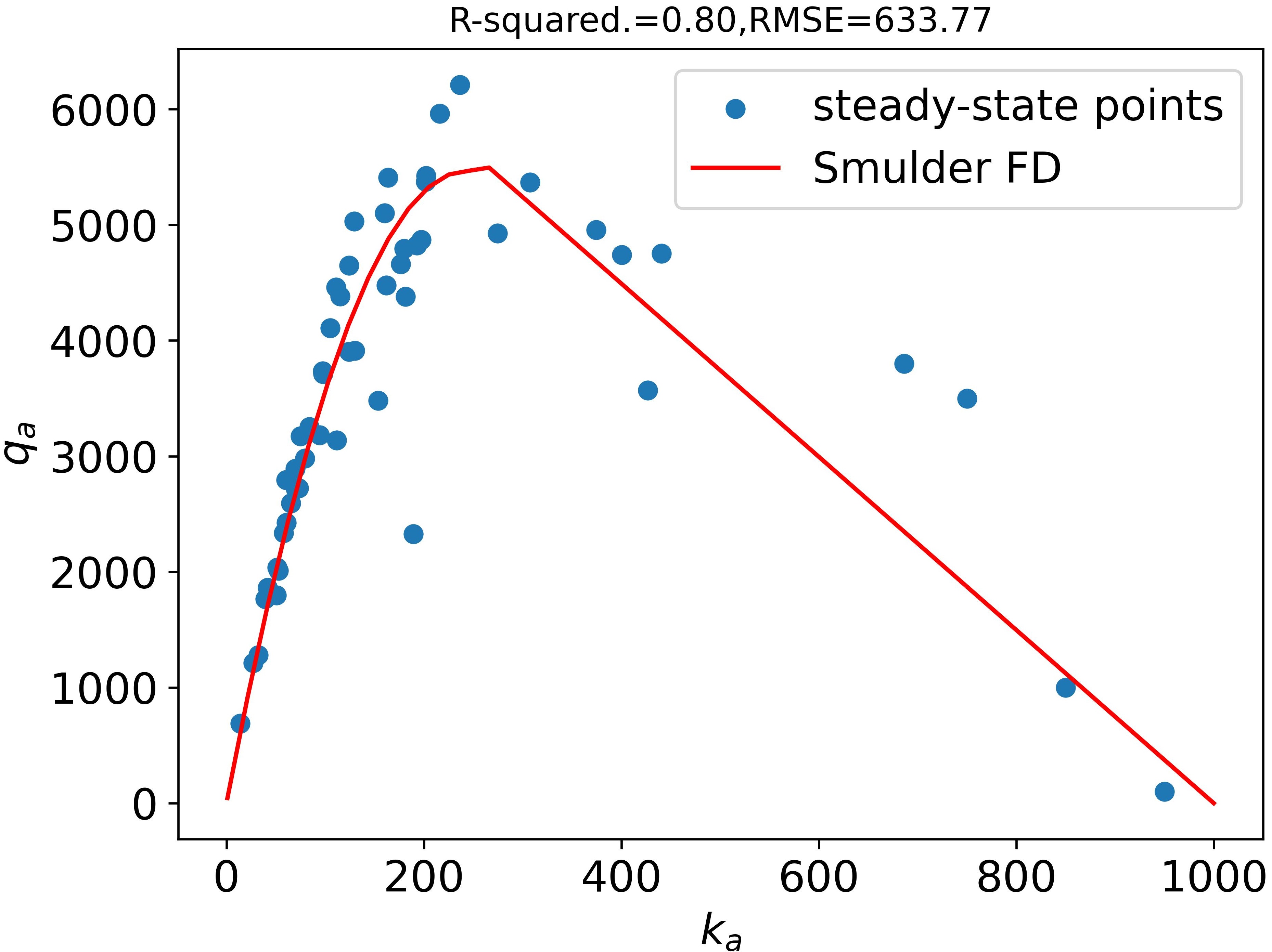}}
\subfigure[Surat-Smulders]{\includegraphics[width=0.3\linewidth]{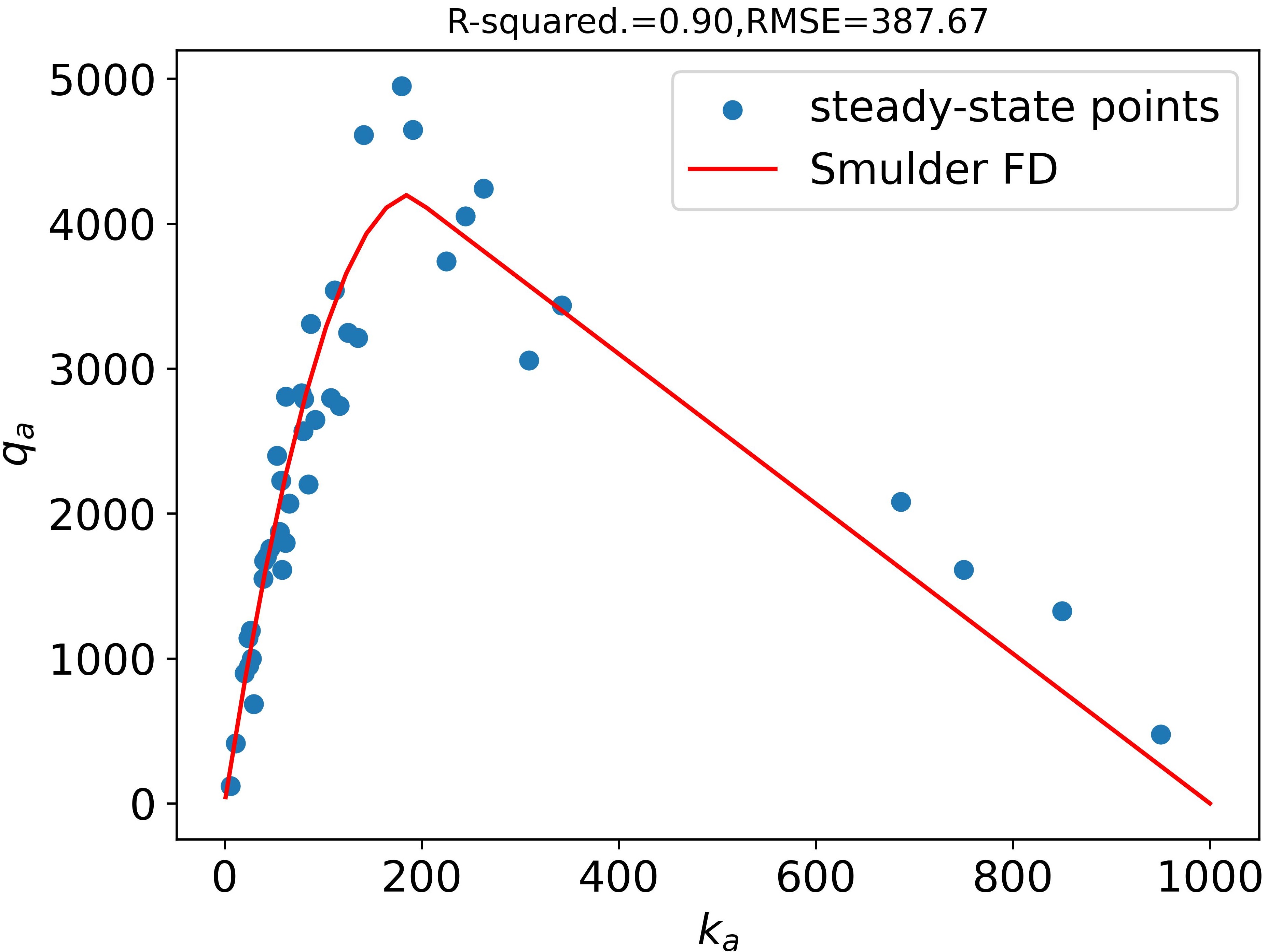}}
\subfigure[Guwahati-Smulders]{\includegraphics[width=0.3\linewidth]{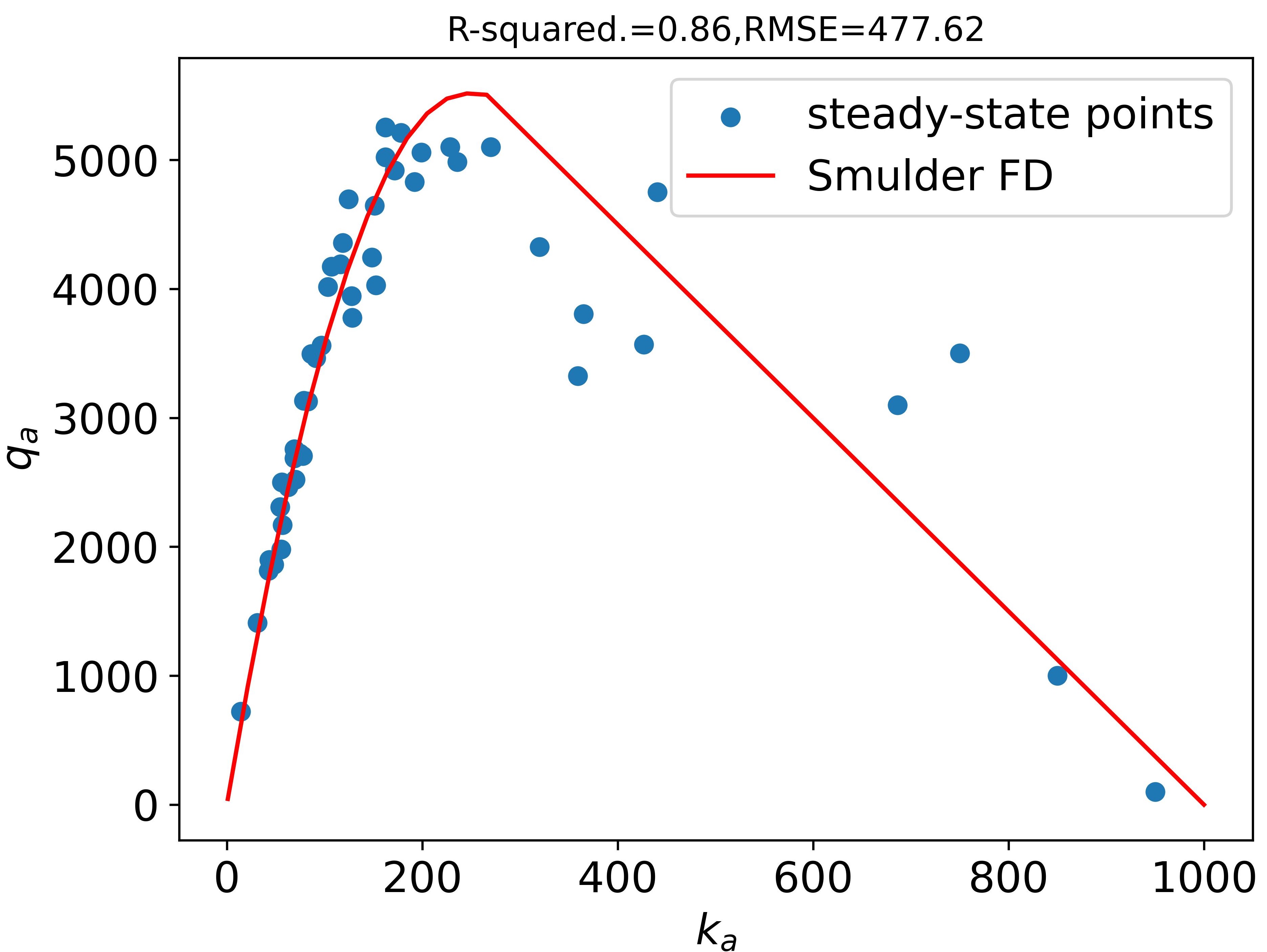}}
\caption{Stream scattered points fitted with Smulders model,  where the left column, middle column, and right column show the Chennai, Surat, and Guwahati data consecutively.}
\label{fig:stream_model}
\end{figure}

 To visualize the relationships among areal flow, areal density, and areal speed of the traffic stream, we plotted a scatter plot of near-stationary empirical data from three distinct locations (see Figure \ref{fig:stream_model}). These empirical data were used to calibrate the fundamental diagram models.

Among the evaluated models, the Smulders model exhibited the best-fitting performance, attributed to its combination of parabolic and linear piecewise functions in the \(q_a-k_a\) plot. This model consistently outperformed others across all three locations, as indicated by higher \(R^2\) values and lower root-mean-square-error (RMSE) values (see Table \ref{tab:model_com_stream}). The scatter plot of steady-state points (Figure \ref{fig:stream_model}) further supports the Smulders model's superior fit. While both Underwood's and Del Castillo's models showed performance levels comparable to Smulders, there were notable differences. Underwood's model estimated a free-flow speed similar to the observed maximum speed but resulted in a significantly higher calibrated critical areal density. Additionally, this model lacks a bounded value for jam density, leading to an inefficient wave speed when compared to observed data. Del Castillo's model, a single-regime model, also aligned well with the empirical data. However, its performance declines due to the assumption of constant speed at low densities. Similarly, Daganzo's double-regime model follows a similar trend to Del Castillo's, differing only in its assumption of constant wave speed in congestion. A nonlinear least squares optimization method was employed to calibrate the models, as detailed in \citet{maiti2024universality}. When comparing the model parameters, the Underwood model demonstrated the highest maximum free-flow speed but failed to accurately estimate critical speed and wave speed. Conversely, the Smulders model provided slightly lower free-flow speed estimates than the Underwood model but delivered the most accurate estimates for both critical speed and wave speed. For calibration purposes, the theoretical value of \(k_{a,jam} = 1000\) was applied. However, obtaining empirical evidence for \(k_{a,jam}\) remains challenging.

\begin{table}[H]
\caption{Comparison of model parameters for the stream with the model fitness score for all data}
\label{tab:model_com_stream}
\begin{tabular}{@{}lllllllllr@{}}
\toprule
\multicolumn{6}{c}{Stream}                                                                                                                    & \multicolumn{2}{c}{$k_a-v_a$}       & \multicolumn{2}{c}{$q_a-k_a$}      \\ \midrule
Location & Model               & $v_{f}$ & $v_{cr}$ & $k_{a,cr}$ & $\omega_a $  & $R^2$         & $RMSE$          & $R^2$         & $RMSE$         \\ \midrule
 & Greenshields                & 38   &   -   &     -  &   -   & 0.72          & 6.59          & 0.09          & 1860         \\
         & Greenberg                   &  -    & 14.5 &  -     &   -   & 0.85          & 4.74          & 0.79          & 652          \\
         & Underwood                   & 50   &  -    & 290   &   -   & 0.91          & 3.82          & 0.80           & 636      \\
Chennai  & Del Castillo                & 40   &    -  &    -   & 6.4  & 0.90           & 3.94          & 0.75          & 719       \\
         & Daganzo                     & 38   &  -    & 135   & 5.9 & 0.87          & 4.56          & 0.68          & 803       \\
         & Smulders                    & 45   & 21   & 255   & 7.5  & \textbf{0.92}         & \textbf{3.77}          & \textbf{0.80}           & \textbf{634}      \\ \hline
         & Greenshields                & 36.5 &    -  &  -     &    -  & 0.7           & 7.27          & 0.07          & 1871         \\
         & Greenberg                   &    -  & 11.9 &  -     &  -    & 0.85          & 5.13          & 0.80           & 531          \\
         & Underwood                   & 46.7 &  -    & 237&   -   & 0.91          & 3.94          & 0.89          & 396        \\
Surat    & Del Castillo                & 39.5 &    - &   -   & 4.3  & 0.85          & 4.97          & 0.79          & 557        \\
         & Daganzo                     & 35.5 &      & 135   & 5    & 0.78          & 6.10           & 0.82          & 514       \\
         & Smulders                    & 43.5 & 21   & 200   & 5.2 & \textbf{0.90}           & \textbf{4.14}          & 0.90           & 388      \\  \hline
         & Greenshields                & 38.7 &   -   &    -   &   -   & 0.78          & 5.74          & 0.08          & 2097         \\
         & Greenberg                   &    -  & 14.5 &  -     & -     & 0.87          & 4.43          & 0.76          & 623          \\
Guwahati & Underwood                   & 49   &   -   & 309   &  -    & 0.95          & 2.81          & 0.71          & 692      \\
         & Del Castillo                & 39.5 &   -   &   -    & 7    & 0.93          & 3.18          & 0.88          & 443       \\
         & Daganzo                     & 35.5 &    -  & 135   & 6.1 & 0.86          & 4.54          & 0.83          & 532        \\
         & \textbf{Smulders}                    & 44.8 & 22.6 & 255   & 7.5  & \textbf{0.95 }        & \textbf{2.64}          & 0.86          & 478      \\ \bottomrule
\end{tabular}
\end{table}

\subsection{Multi-class fundamental diagrams with areal variables}
Now, we analyze the multi-class fundamental diagrams using the proposed areal variables to gain a better understanding of the class-specific vehicle dynamics and their interactions. The existing multiclass models estimate class-specific speed and flow as functions of the total density. Similarly, this section models class-specific speed as a function of the total areal density. For the multi-class speed function, we calibrated all the functional forms (see in Table \ref{tab:speed_density_models}) of the speed function for each vehicle class (cars, TWs, HVs) as a function of total areal density. The calibrated parameters, along with their fitting scores, are summarized in  Table \ref{tab:model_com_tw} for two-wheelers, Table \ref{tab:model_com_car} for cars, Table \ref{tab:model_com_hvs} for heavy vehicles. Among the  models assessed, the Smulder model consistently exhibited superior performance across all vehicle classes and data locations, as evidenced by its lower RMSE and higher $R^2$ values. Therefore, we selected the Smulder function for the multiclass derivation and further numerical analysis.    
The generic class-specific flow and speed function for multi-class FDs is defined in \eqref{eq:mc_exp}.
\begin{equation}\label{eq:mc_exp}
    q_a^i=k_a^iv^i, \hspace{1cm}  v^i=v^i(k_a), \hspace{1cm} k_a=\sum_{\forall i}  k_a^i
\end{equation}

The proposed multi-class areal density-areal speed fundamental relationships for Smulders model are shown in \eqref{eq:str_fd_smulder1}, \eqref{eq:str_wave1}. Class-specific speed $(v^i(k_a))$ is a decreasing function of total areal density ($k_a$), ranging between \textit{i}-th class vehicle's maximum observed speed $(v_{f}^i)$ to critical speed $(v_{cr}^i)$ below \textit{i}-th class critical areal density $(k_{a,cr}^i)$. Beyond $k_{a,cr}^i$, $v^i(k_a)$ is a function of class-specific wave speed $(\omega_i)$ and total areal density $(k_a)$. The variables $v_{f}^i, v_{cr}^i, k_{a,cr}^i, k_{a,jam}$, and  $ \omega_a^i$ are the model parameter and do not depend on traffic state. The multi-class FDs for the Smulders model are demonstrated in \eqref{eq:str_fd_smulder1}, \eqref{eq:str_wave1}. These class-specific FDs also satisfy the general principle of the model presented in \eqref{eq:mc_exp}. This study relaxes the requirement of class-independent wave speed as proposed by \citet{van2014extension} and proposed class-specific wave speed in congestion, $ \omega_a^i$.
\newline Smulders class-specific FDs:
\begin{equation} \label{eq:str_fd_smulder1}
    v^i(k_a)= 
\begin{cases}
    v_{max}^i-\frac{(v_{max}^i-v_{crit}^i)k_a}{k_{a,crit}^i},& \text{if } 0\leq k_a \leq k_{a,crit}^i\\
    \omega_a^i\left(\frac{k_{a,jam}}{k_a}-1\right),              & \text{otherwise}
\end{cases}
\end{equation}
Class-specific wave speed:
\begin{equation}\label{eq:str_wave1}
   \omega_a^i =\frac{v_{crit}^i k_{a,crit}^i}{k_{a,jam}-k_{a,crit}^i}
\end{equation}
Now, we evaluated the fitness of the proposed multiclass model against empirical data to validate the model parameters. Table \ref{tab:model_com_smulders} presents the class-specific model parameters estimated under the assumption that class-specific variables adhere to the Smulders fundamental relationships, along with their goodness-of-fit measures when compared to the empirical dataset. The model parameters were calibrated using speed-density relationships, and these optimized parameters were subsequently used to estimate flow-density models. This calibration approach accounts for the higher \( R^2 \) values observed for the speed-density relationship \((k_a-v_a^i)\) compared to the flow-density relationship \((k_a-q_a^i)\).

Table \ref{tab:model_com_smulders} provides a detailed comparison of parameters such as the maximum observed speed \((v_{f}^i)\), critical speed \((v_{cr}^i)\), critical areal density \((k_{a,cr}^i)\), and wave speed \((\omega_a^i)\) across different vehicle classes. Notably, the maximum observed speed \((v_{f}^i)\) is almost identical for each vehicle class across all locations. Interestingly, heavy vehicles exhibit higher critical density and wave speeds compared to two-wheelers and cars across all locations. The order of critical density and wave speed follows the pattern: 

\[
k_{a,cr}^{\text{HVs}} > k_{a,cr}^{\text{cars}} > k_{a,cr}^{\text{TWs}}
\]

Similarly, for wave speed:

\[
\omega_{a}^{\text{HVs}} > \omega_{a}^{\text{cars}} > \omega_{a}^{\text{TWs}}
\]

A comprehensive comparison of class-specific fundamental diagrams for various FD models across all study locations is provided in Appendix Tables \ref{tab:model_com_tw}, \ref{tab:model_com_car}, and \ref{tab:model_com_hvs}.
The results of this study support the findings of \citet{VanWageningen-Kessels2018}, demonstrating that the vehicle class with the largest vehicles, rather than the longest vehicles, exhibited the highest wave speed.  This observation indicates that the new multi-class model with areal variables aligns with traditional multi-class models in terms of wave propagation and critical density characteristics in the traffic stream.

\par The results of analyzing stream and class-specific FDs suggest that newly introduced areal variables are meaningful and valuable in portraying real-world equilibrium conditions in mixed traffic settings. Furthermore, we observed that the individual class fundamental diagrams exhibit minimal variation within the stream. Also, the comparison between the existing stream models and the class-specific models indicates that the stream representation yields optimal values for both the coefficient of determination ($R^2$) and root-mean-square-error (RMSE) for both the $k_a$-$v$ and $k_a$-$q_a$ fundamental diagrams. The study's findings suggest that a stream flow model could serve as an appropriate representation of heterogeneous urban corridor traffic. Moving forward, we delve into examining the conventional solution approaches of the kinematic wave model within the framework of the proposed continuum model featuring areal fundamental diagrams.

\begin{figure}[H]
\centering
\subfigure[Chennai-Smulders]{\includegraphics[width=0.3\linewidth]{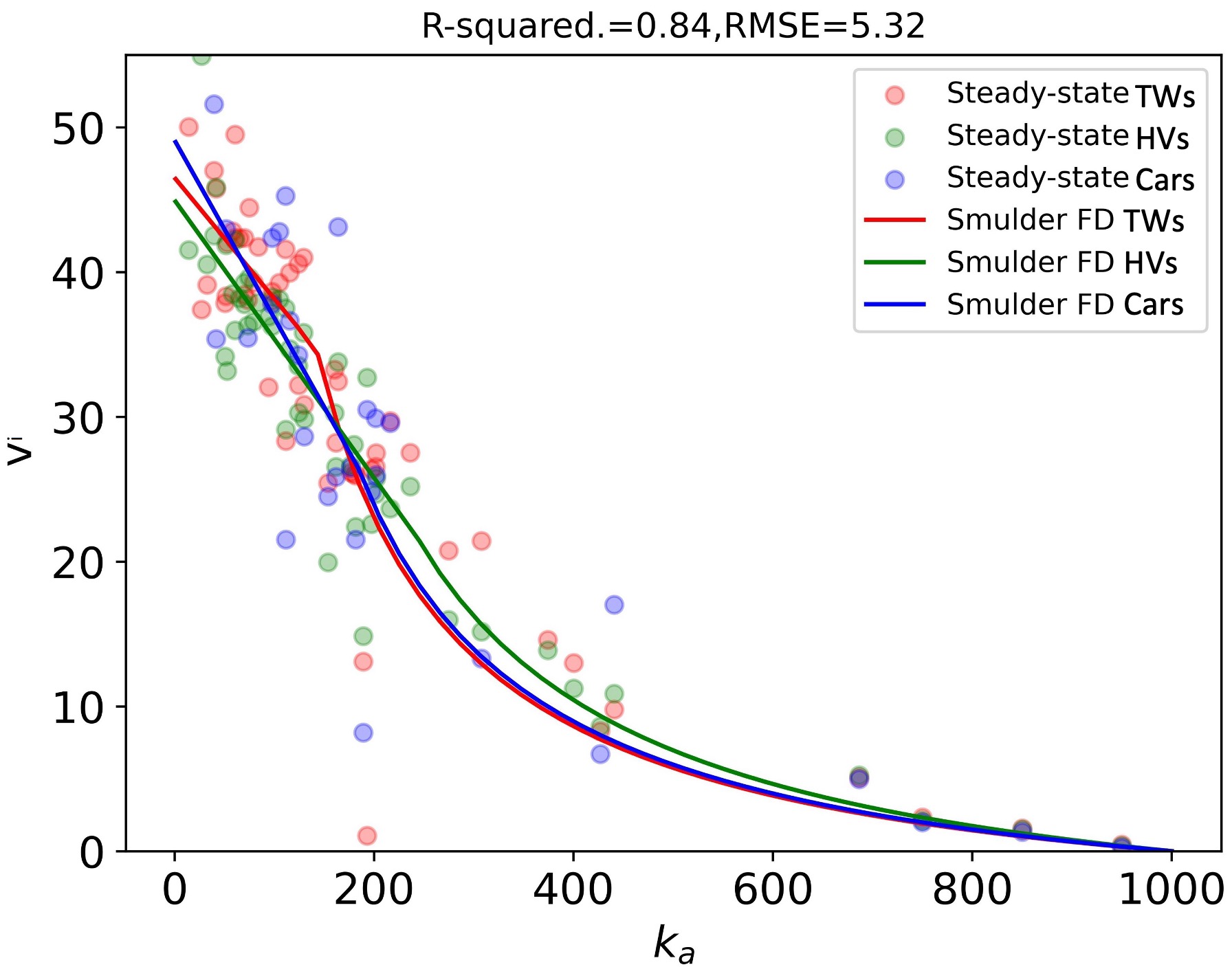}}
\subfigure[Surat-Smulders]{\includegraphics[width=0.3\linewidth]{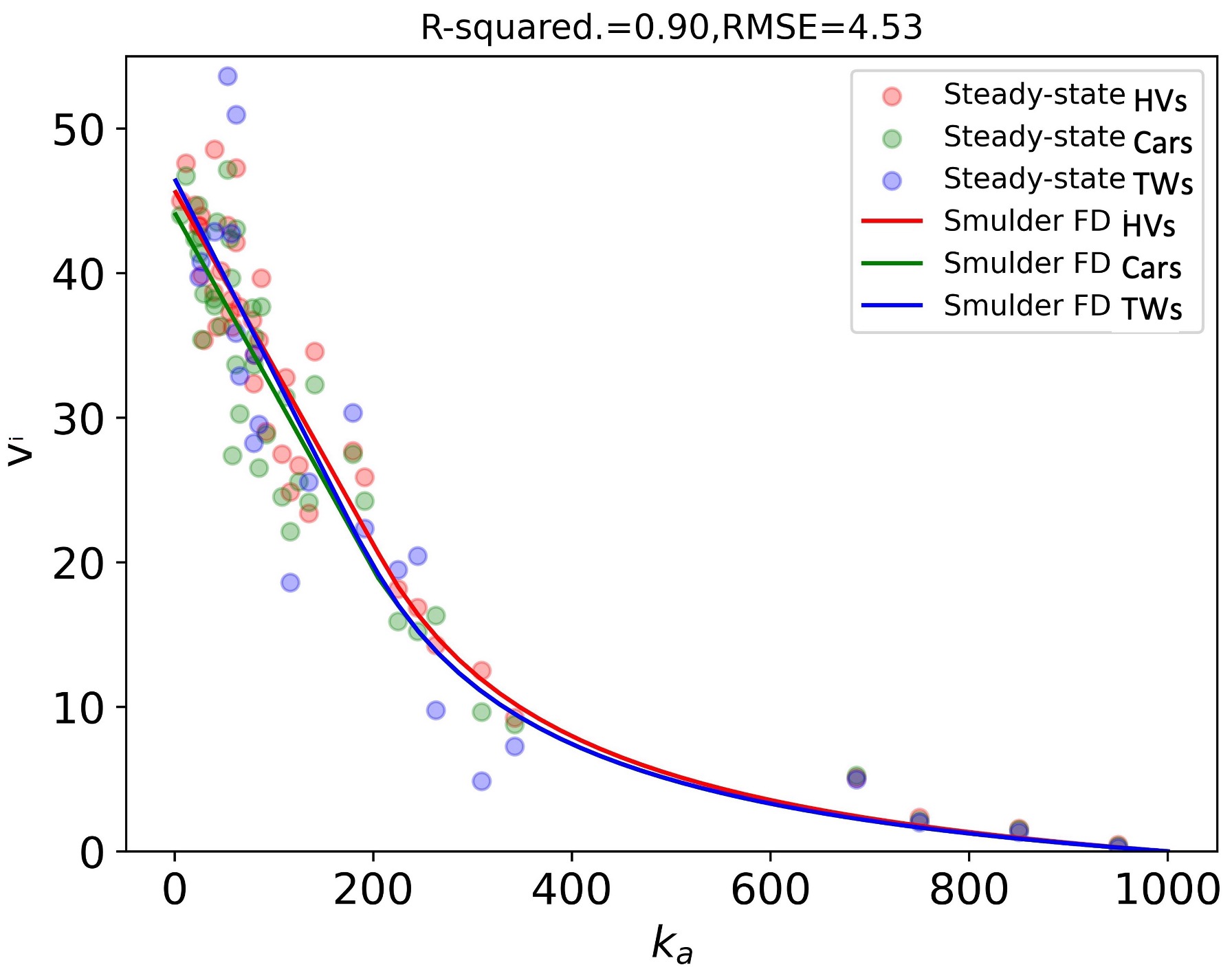}}
\subfigure[Guwahati-Smulders]{\includegraphics[width=0.3\linewidth]{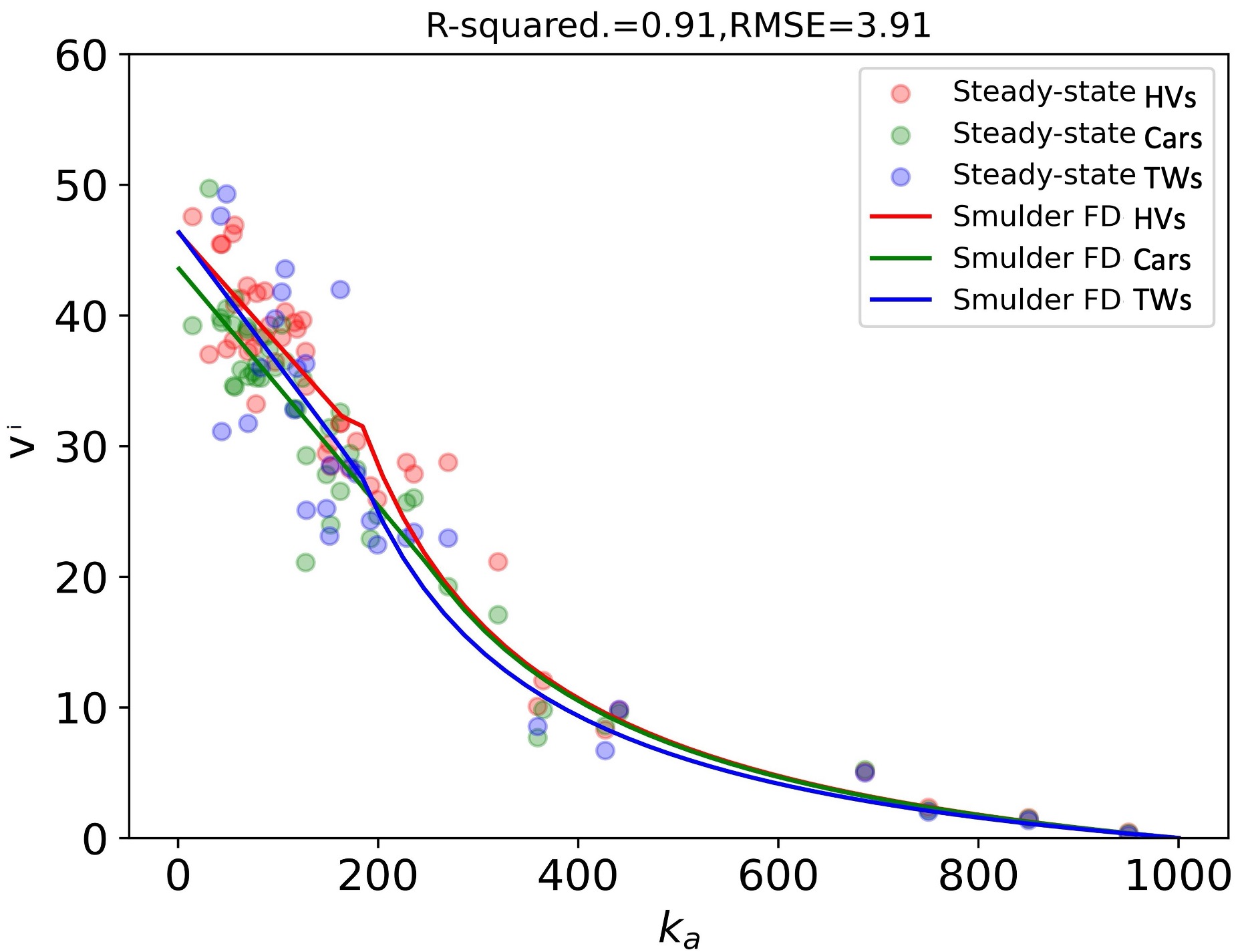}}
\subfigure[Chennai-Smulders]{\includegraphics[width=0.3\linewidth]{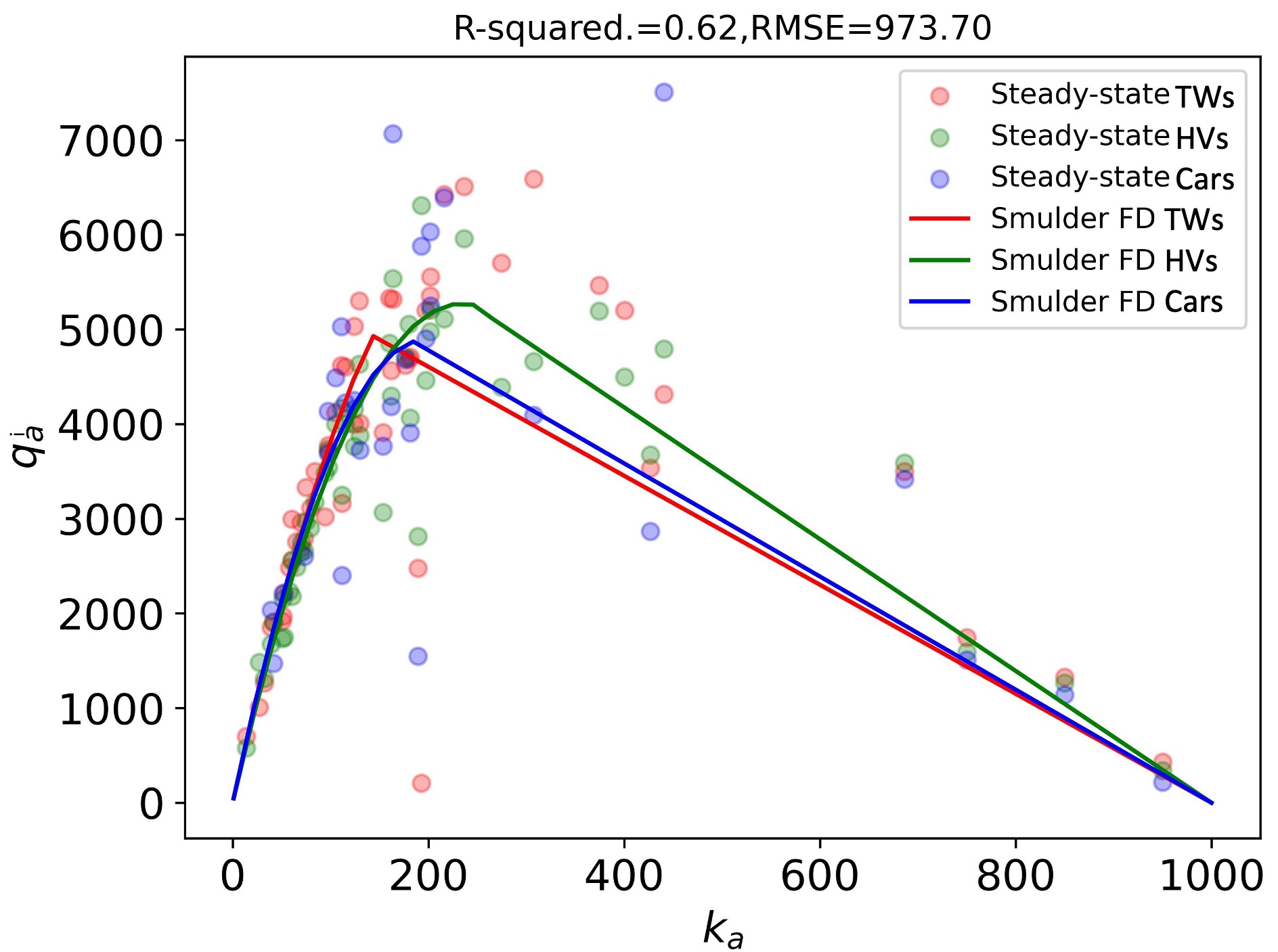}}
\subfigure[Surat-Smulders]{\includegraphics[width=0.3\linewidth]{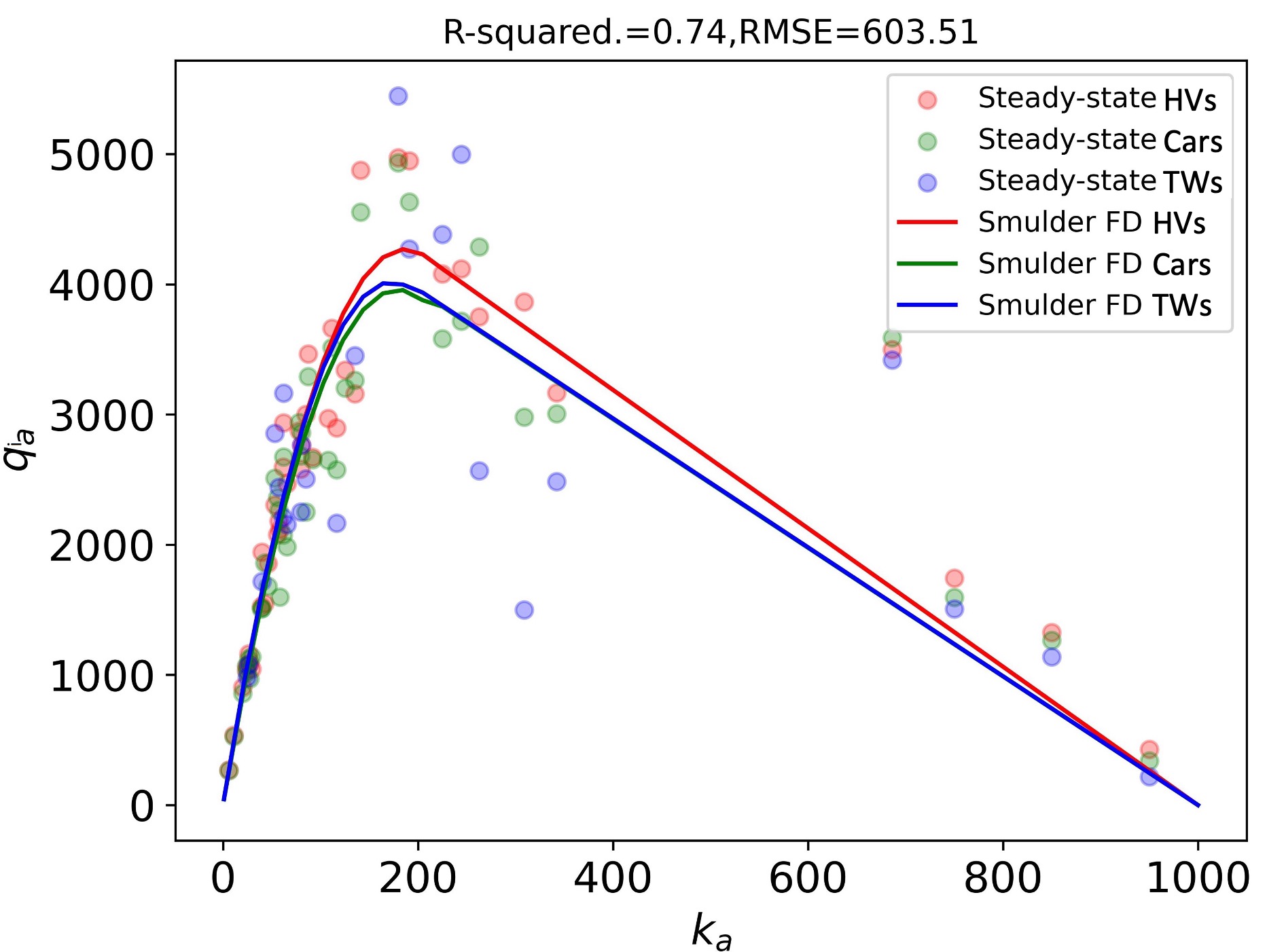}}
\subfigure[Guwahati-Smulders]{\includegraphics[width=0.3\linewidth]{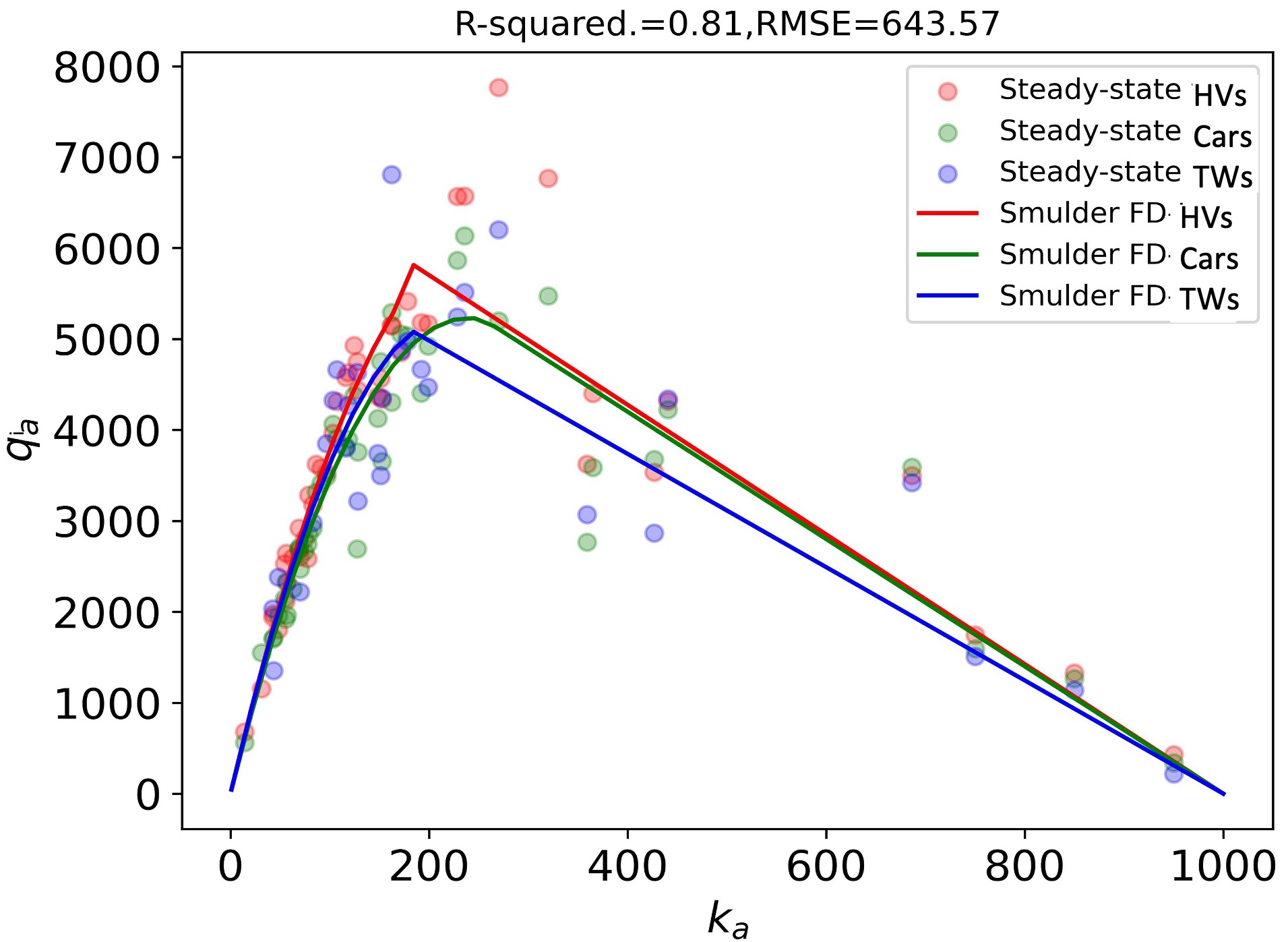}}
\caption{Scatter steady-state points fitted with the Smulders and proposed fundamental relationships for multiple classes, where the left, middle, and right columns show the Chennai, Surat, and Guwahati data consecutively. }
\label{fig:multiclass_model}
\end{figure}

\begin{table}[H]
\caption{Comparison of Smulders model parameters for two-wheelers (TWs), cars, and heavy vehicles (HVs) with the model fitness score for all data}
\label{tab:model_com_smulders}
\begin{tabular}{@{}lllllllll@{}}
\toprule
\multicolumn{5}{c}{Smulders model}                                     & \multicolumn{2}{c}{$k_a-v^i_a$}       & \multicolumn{2}{c}{$q^i_a-k_a$}      \\ \midrule
Location            & $v^i_{f}$ & $v^i_{cr}$ & $k^i_{a,cr}$ & $\omega^i_a$   & $R^2$         & $RMSE$          & $R^2$         & $RMSE$         \\ \midrule
 & & & & TWs & & & & \\ \midrule
      Chennai  & 49.5 & 29.00 & 170.00   & 5.94 & 0.84          & 5.32          & 0.62  & 973          \\ 
      Surat  & 48.40 & 19.50 & 193.00   & 4.66 & 0.90          & 4.53          & 0.74          & 603          \\ 
     Guwahati & 48.80 & 29.40 & 170.00   & 6.02 & 0.90          & 3.93          & 0.80          & 647          \\ \midrule
 & & & & Cars & & & & \\ \midrule
 Chennai & 49.00 & 25.00  & 200.00   & 5.25 & 0.82          & 6.50          & 0.61          & 980          \\ 
 Surat & 48.80 & 19.70  & 195.00   & 4.77 & 0.89          & 4.80          & 0.73          & 612          \\
Guwahati & 48.80 & 25.60  & 195.00   & 6.20 & 0.91          & 3.92          & 0.81          & 643          \\ \midrule
 & & & & HVs & & & & \\ \midrule
  Chennai& 49.00 & 23.00 & 250.00   & 7.67 & 0.84          & 5.40          & 0.62          & 975          \\ 
  Surat & 48.20 & 20.00 & 210.00   & 5.32  & 0.89          & 4.60          & 0.72          & 604          \\
  Guwahati& 48.65   & 30.00   & 200.00   & 7.50  & 0.89          & 3.96          & 0.79          & 645          \\ \bottomrule
\end{tabular}
\end{table}

\section{Numerical Solution}
This section proposes a multi-class cell transmission model to demonstrate the applicability of the proposed areal continuum model to generate mixed traffic scenarios, including platoon dispersion and overtaking phenomena. Also, we illustrated the applicability of the existing numerical methods for the hyperbolic partial differential equation (PDE) to the proposed areal model.
\subsection{Shock-Wave Analysis}
\begin{figure}[h]
\centering  
\subfigure[]{\includegraphics[width=0.4\linewidth]{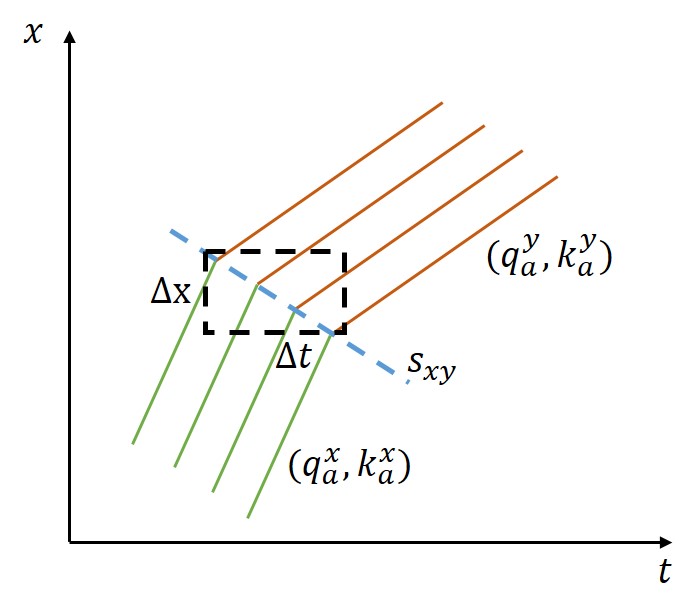}}
\subfigure[]{\includegraphics[width=0.4\linewidth]{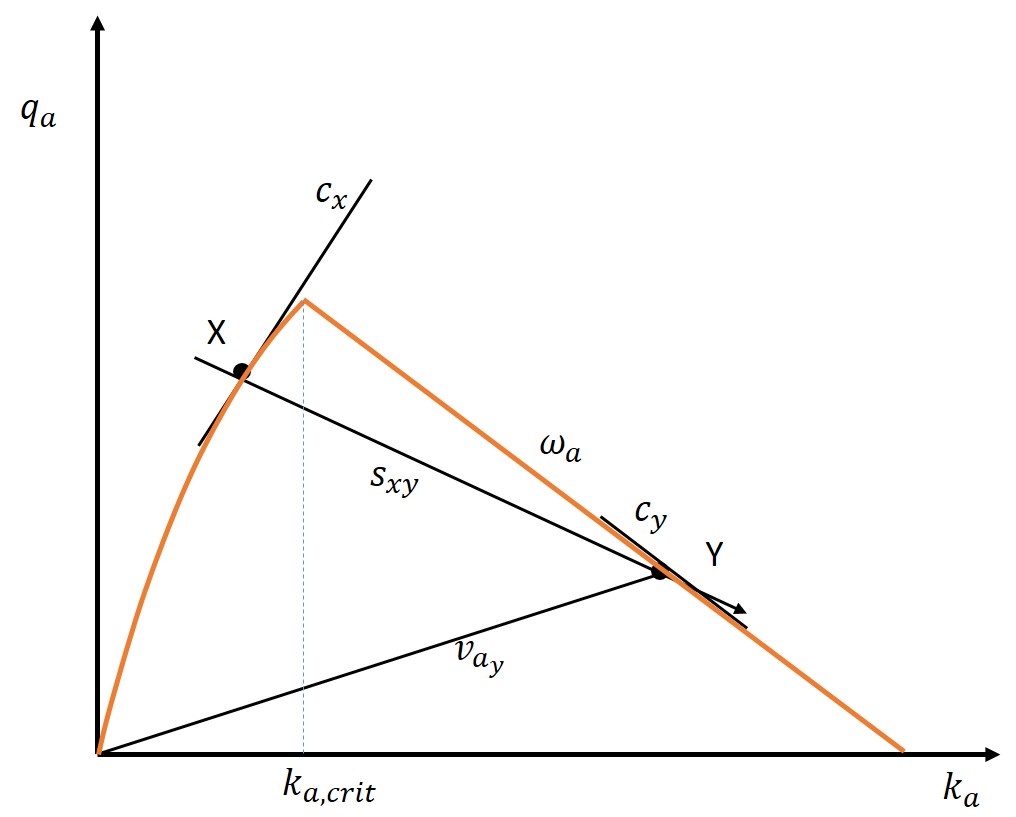}}
\caption{Application of $q_a$-$k_a$ curve. The slope of the radius vector at point y gives the average AS $(v_{Ay})$ of vehicles; the slope of the tangent at x, y gives wave speed $c_x,c_y$; line joining two traffic states x, y gives the local speed near a shock-wave }\label{fig:wave}
\end{figure}

The speed of a shock wave, in the traditional sense, is defined as the proportion of the relative flow to the relative density between two interacting traffic states. This relationship can be visually represented by traffic state x meeting with a more dense state y, producing a backward-moving shock with a speed of $s_{xy}$ (in Figure \ref{fig:wave}).  
\newline Theorem: The shock speed is the slope of the line joining the upstream and downstream traffic states in the areal ($q_a$-$k_a$) fundamental diagrams.   
\par 
Based on the principle of vehicle area conservation, the total vehicle area entering the box (see Figure \ref{fig:wave} (a)) is equal to the total vehicle area getting out from the box. 
\begin{equation}\label{eq:shock}
    \begin{split}
        Veh\_area\_in &=Veh\_area\_out \\
        k_a^x \Delta x w+q_a^x \Delta t w&=k_a^y \Delta x w+q_a^y \Delta t w\\
        -\frac{\Delta x}{\Delta t}&=\frac{q_a^x-q_a^y}{k_a^x-k_a^y}\\
        s_{xy} &=\frac{q_a^x-q_a^y}{k_a^x-k_a^y}
    \end{split}
\end{equation}

Where $k_a^x$, $k_a^y$ is the areal density at upstream and downstream traffic states, and $q_a^x$, $q_a^y$ is the corresponding areal flow. Eq. \eqref{eq:shock} reveals the magnitude of shock speed is equal to the proportion of relative areal flow to the relative areal density of two interacting traffic states.

At a point on the curve, the slope of the tangent to the curve reveals the speed of waves, $c$ ($c_{x/y}$ wave speed at x or y point in Figure \ref{fig:wave} (b)), illustrate continuous changes in flow through the stream. In the free flow regime, the wave speed transitions from the maximum free flow speed $(v_{max})$ at low densities to a negative congestion wave speed $(-\omega_a)$ as the traffic becomes congested. The expression for wave speed, as given by \eqref{eq:wave-eqn}, is a function of $k_a$ in the free flow regime and remains constant as $\omega_a$ in the congested regime.
\begin{equation} \label{eq:wave-eqn}
    c(k_a)= \frac{\partial q_a}{\partial k_a}=
\begin{cases}
    v_{f}-2k_a\frac{v_{f}-v_{cr}}{k_{a,cr}},
    \text{if } k_a \leq k_{a,cr}\\
    -\omega_a,
    \text{otherwise}
\end{cases}
\end{equation}
Areal speed, $v_a=\frac{q_a}{k_a}$, represents the slope of the radius vector from the origin, which equals the space-mean speed at any point on the curve.
\subsection{Analytical Solution: Method of Characteristics}
\begin{figure}[h]
  \centering
  \includegraphics[width=0.8\textwidth]{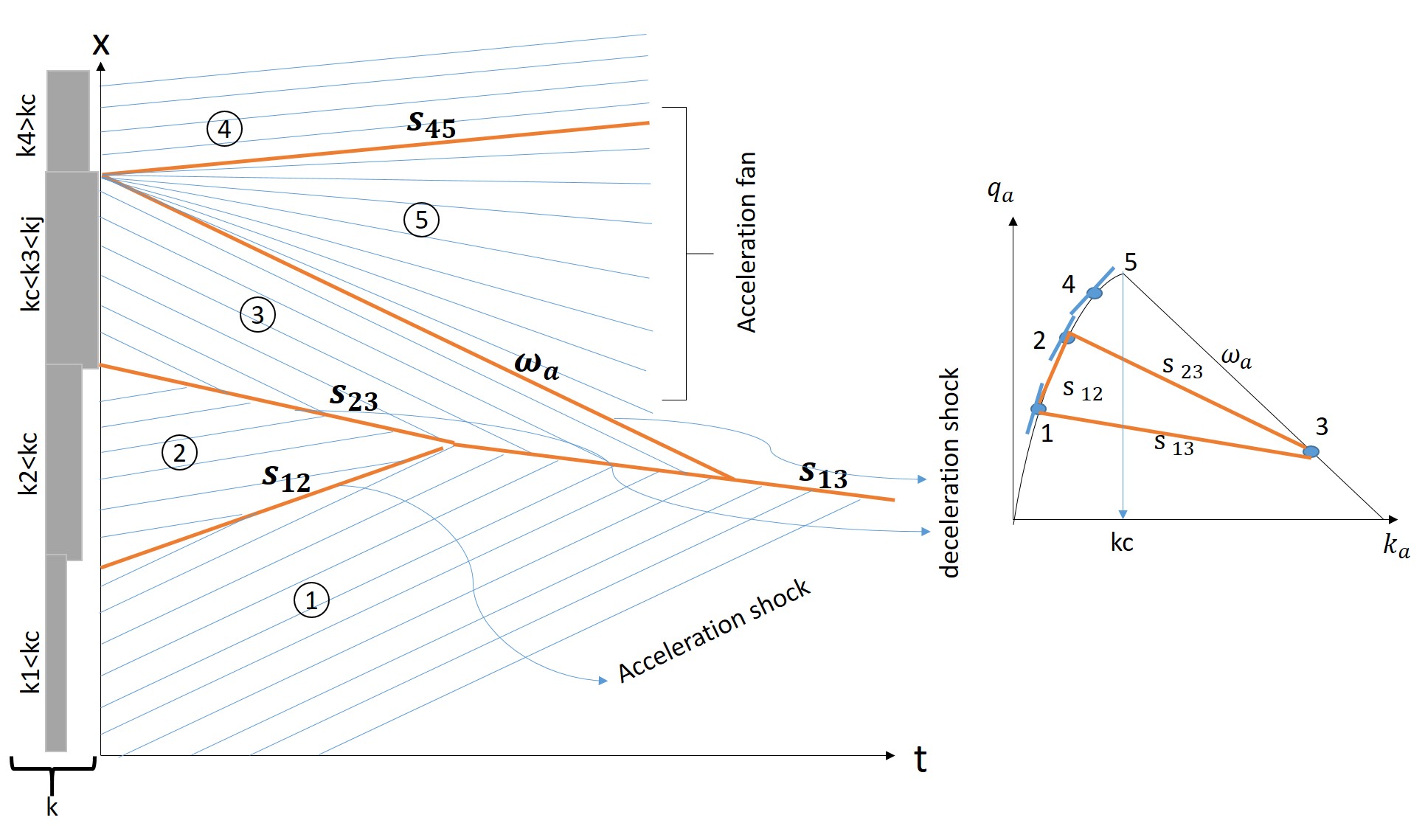}
  \caption{Characteristic lines and shock speeds for a initial traffic conditions}\label{fig:moc_exp}
\end{figure}

The  areal continuum model \eqref{eq:final_cont} can be described as \eqref{eq:analy2}.
Now substitute \eqref{eq:wave-eqn} into the continuity equation \eqref{eq:final_cont}, we get
\begin{equation}\label{eq:analy2}
    \frac{\partial k_{a}}{ \partial t}+\frac{dq_a}{dk_a}\frac{\partial k_{a}}{\partial x }=0
\end{equation}

\begin{equation}\label{eq:wave}
  c(k_a)=\frac{dq_a}{dk_a}
\end{equation}
Note that the general form of the \eqref{eq:analy2} after applying the substitute \eqref{eq:wave} as follows:
\begin{equation}
  \frac{\partial k_{a}}{ \partial t}+c(k_a)\frac{\partial k_{a}}{\partial x }=0  
\end{equation}
This is a quasilinear PDE, just like the LWR model. Therefore, the method of characteristics (MoC) can be effectively applied to the areal continuum model to derive analytical solutions.
An example of MoC is shown in Figure \ref{fig:moc_exp}. Characteristic lines are drawn from $t=0$ across the length of the stretch $(x \in [0, L])$, based on the initial traffic conditions and utilizing Smulder's stream flow-density relationships. he figure depicts examples of  acceleration shocks $(s_{12},s_{45})$, and deceleration shocks $(s_{13}, s_{23}, \omega_a)$, which occur during transitions between different traffic states. These transitions are clearly delineated in Figure \ref{fig:moc_exp}, highlighting the dynamics of traffic state propagation as captured by the MoC using the proposed model.

\subsection{Numerical solution: kinematic wave (KW)-based approach }
The newly developed continuum model belongs to the hyperbolic PDE class, like the traditional LWR models. Hence, the traditional numerical schemes can be used to simulate heterogeneous traffic systems.
\newline The analytical results presented in section 4.2 provide insights into the formation and dissipation of queues and congestion over time and space. However, to simulate complex phenomena, like long highways with complex geometry and traffic conditions with multiple sinks and sources, a numerical methodology is needed, allowing for the inclusion of these complexities encountered in real traffic systems.

\subsubsection{Multiclass cell transmission model}
\begin{figure}[h]
  \centering
  \includegraphics[width=0.5\textwidth]{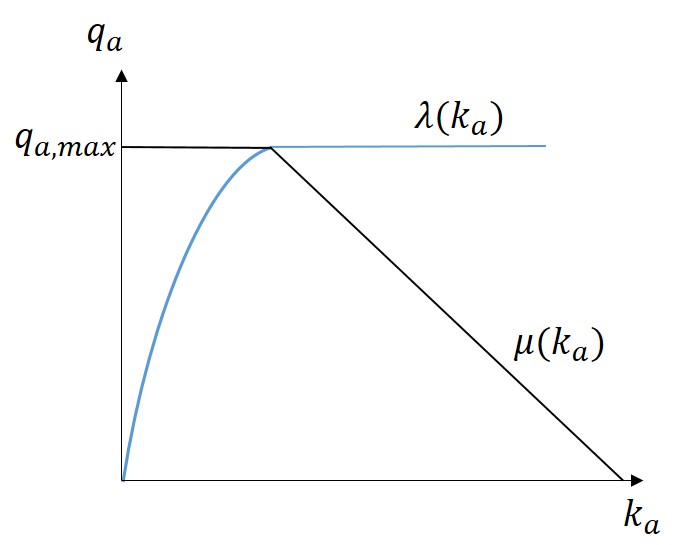}
  \caption{Demand-supply functions for the Godunov's scheme as per the proposed fundamental diagram}\label{fig:godu}
\end{figure}

In the cell transmission model (CTM), \citet{daganzo1994cell} proposed supply $(\mu(k_a))$ and demand function $(\lambda(k_a))$ to estimate the average flow $(\Phi(.))$ at any cell by solving the Riemann problem. 
\begin{equation}
    \Phi_j(k_{a,(j-1)},k_{a,j})=min\{ \lambda(k_{a,(j-1)}),\mu(k_{a,j}), q_{a,max}\}
\end{equation}
The demand-supply function can be defined by the fundamental diagrams. Figure \ref{fig:godu} represents the demand-supply function for this study. Since we used the Smulder FDs in this study, the demand-supply function can be defined as \eqref{eq:dem_sup}.
\begin{equation}\label{eq:dem_sup}
   \Phi_j(k_{a,(j-1)},k_{a,j})=min\{ v(k_{a(j-1)})k_{a(j-1)}),(K_{a,jam}-k_{a,j})\omega_a, q_{a,max}\} 
\end{equation}
Next, we discuss a simple merging model  \citet{daganzo1995cell} to formulate the cell transmission model for multiclass traffic conditions. We consider two classes of vehicles for the simplicity of demonstration: cars and HVs. 
The boundary equations required are the advancing flows, $q_a^{car}(t)$ and $q_a^{HVs}(t)$, in terms of the maximum flow that the sending cell can discharge $(\lambda^{HVs},\lambda^{car} )$  and the immediate downstream capacity of receiving cell, $ \mu^{max}$. Hence, the flows must satisfy
\begin{subequations}\label{eq:mc_ctm_1}
    \begin{align}
      q_a^{car}(t) \leq \lambda_j^{car}; & q_a^{HV}(t) \leq \lambda_j^{HVs} \hspace{1cm} and\\
      q_a^{car}(t)+ q_a^{HV}(t)  &\leq \mu_{j+1}^{max}
    \end{align}
  \end{subequations}

As long as the supply of vehicles from both class $\lambda_j^{car}$ and $\lambda_j^{HV}$ is not exhausted, we will assume that fraction, $p^{car}$ of the vehicles comes from car and the remainder $p^{HV}$ from heavy vehicles, where $p^{car}+p^{HV}=1$.
\begin{figure}[h]
  \centering
  \includegraphics[width=0.45\textwidth]{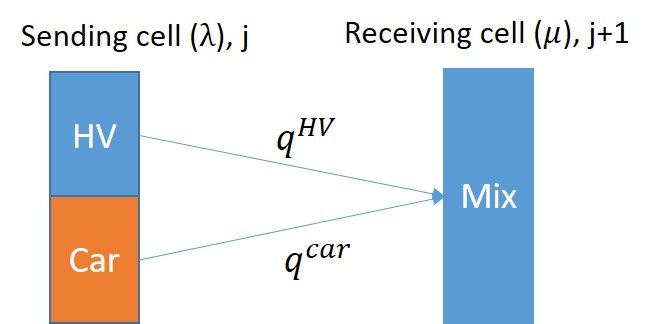}
  \caption{Representation of sending cell and receiving cell of multiclass cell transmission model}\label{fig:mc_ctm}
\end{figure}
Thus, the areal flow for car and HV can be represented by \eqref{eq:mc_ctm3}:
\begin{subequations}\label{eq:mc_ctm3}
    \begin{align}
        q_{a,j}^{car}(t) &= min \{\lambda_j^{car}, (\mu^{mix}_{j+1} - \lambda_j^{HV}), p_j^{car} \mu^{mix}_{j+1}, q_{a,max} \} \hspace{2mm} and\\
      q_{a,j}^{HV}(t) &= min \{\lambda_l^{HV}, (\mu^{mix}_{j+1} - \lambda_j^{car}), p_j^{HV} \mu^{mix}_{j+1}, q_{a,max} \}, \hspace{2mm} if  \hspace{1mm} \mu^{mix}_{j+1} < \lambda_j^{car}+\lambda_j^{HV}
    \end{align}
  \end{subequations}
Further, we assumed the fraction of cars and HVs advance to the $(j+1)$  in \eqref{eq:mc_ctm3} is the function of density proportions of those vehicles at $j$ cell. A similar kind of assumption for density proportion in ramp metering can be found in \citet{gomes2006optimal}.
\begin{subequations}\label{eq:mc_ctm4}
    \begin{align}
        p^{car}_j(t) =\frac{k^{car}_{a,j}(t)}{k^{car}_{a,j}(t)+k^{HV}_{a,j}(t)} \\
      p^{HV}_j(t) =\frac{k^{HV}_{a,j}(t)}{k^{car}_{a,j}(t)+k^{HV}_{a,j}(t)}
    \end{align}
  \end{subequations}
Equation \eqref{eq:mc_ctm3} described the flow update in the cell, while the associated sending and receiving functions are shown in \eqref{eq:dem_sup}. Once the flows are updated to the downstream cell, the subsequent step involves the estimation of density propagation in the next time step using the following equation \eqref{eq:mc_ctm5}. Hence, the multiclass traffic density propagation is achieved by utilizing the merged model proposed by \citet{daganzo1995cell}.
\begin{subequations}\label{eq:mc_ctm5}
    \begin{align}
  k_{a,j}^{car}(t+1)=k_{a,j}^{car}(t)+\frac{\Delta t}{l_j}(q_{a,j+1}^{car}(t)-q_{a,j}^{car}(t))\\
  k_{a,j}^{HV}(t+1)=k_{a,j}^{HV}(t)+\frac{\Delta t}{l_j}(q_{a,j+1}^{HV}(t)-q_{a,j}^{HV}(t))\\
   k_{a,j}(t+1)=k_{a,j}^{car}(t+1)+ k_{a,j}^{HV}(t+1)
\end{align}
  \end{subequations}

\begin{figure}[H]
\centering  
\subfigure[Two vehicle class platoons  $(k^{car}_a=150, k^{car}_a=150)$ at a particular location of a road, create congestion ]{\includegraphics[width=0.3\linewidth]{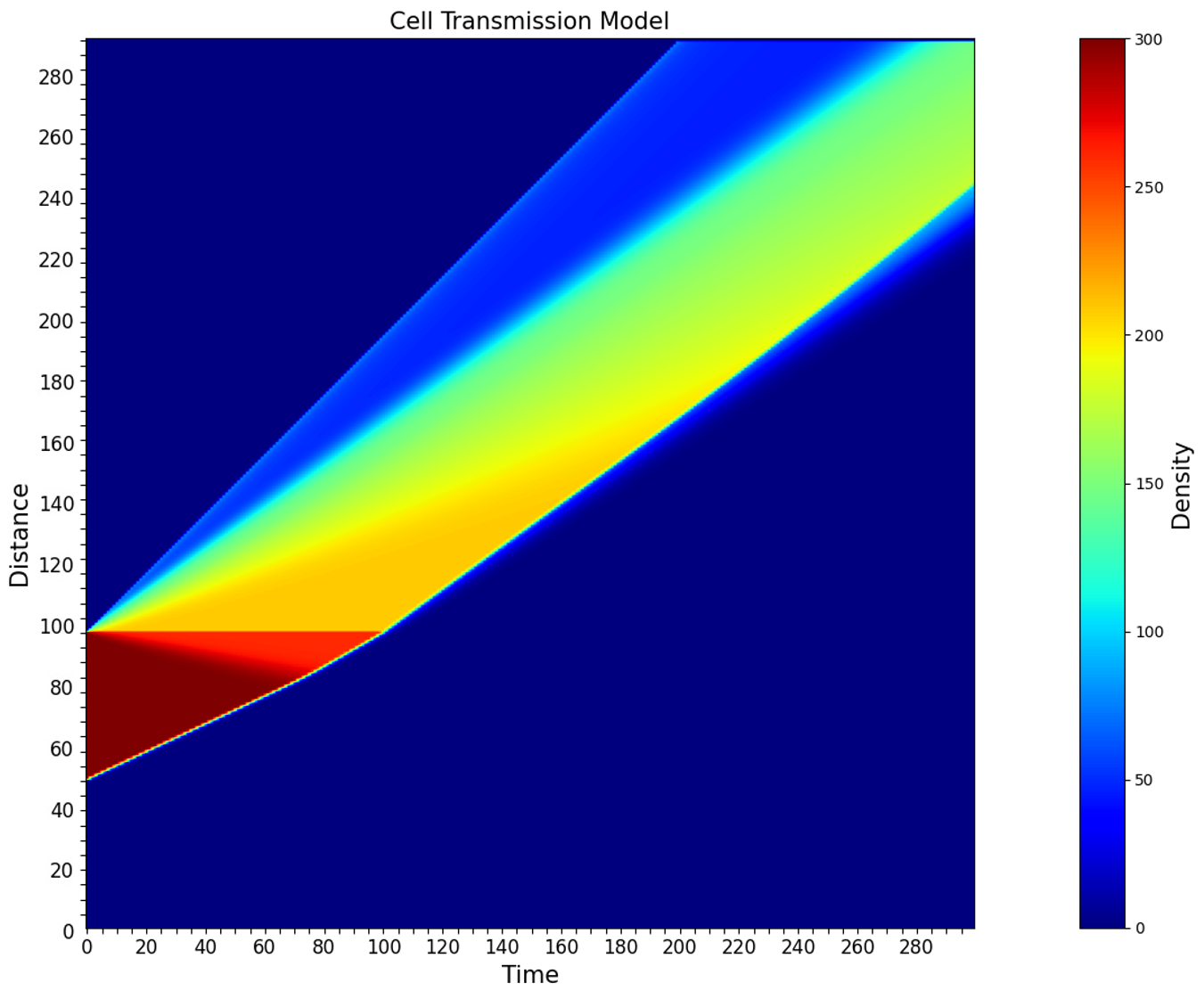}}
\subfigure[$t=0$]{\includegraphics[width=0.22\linewidth]{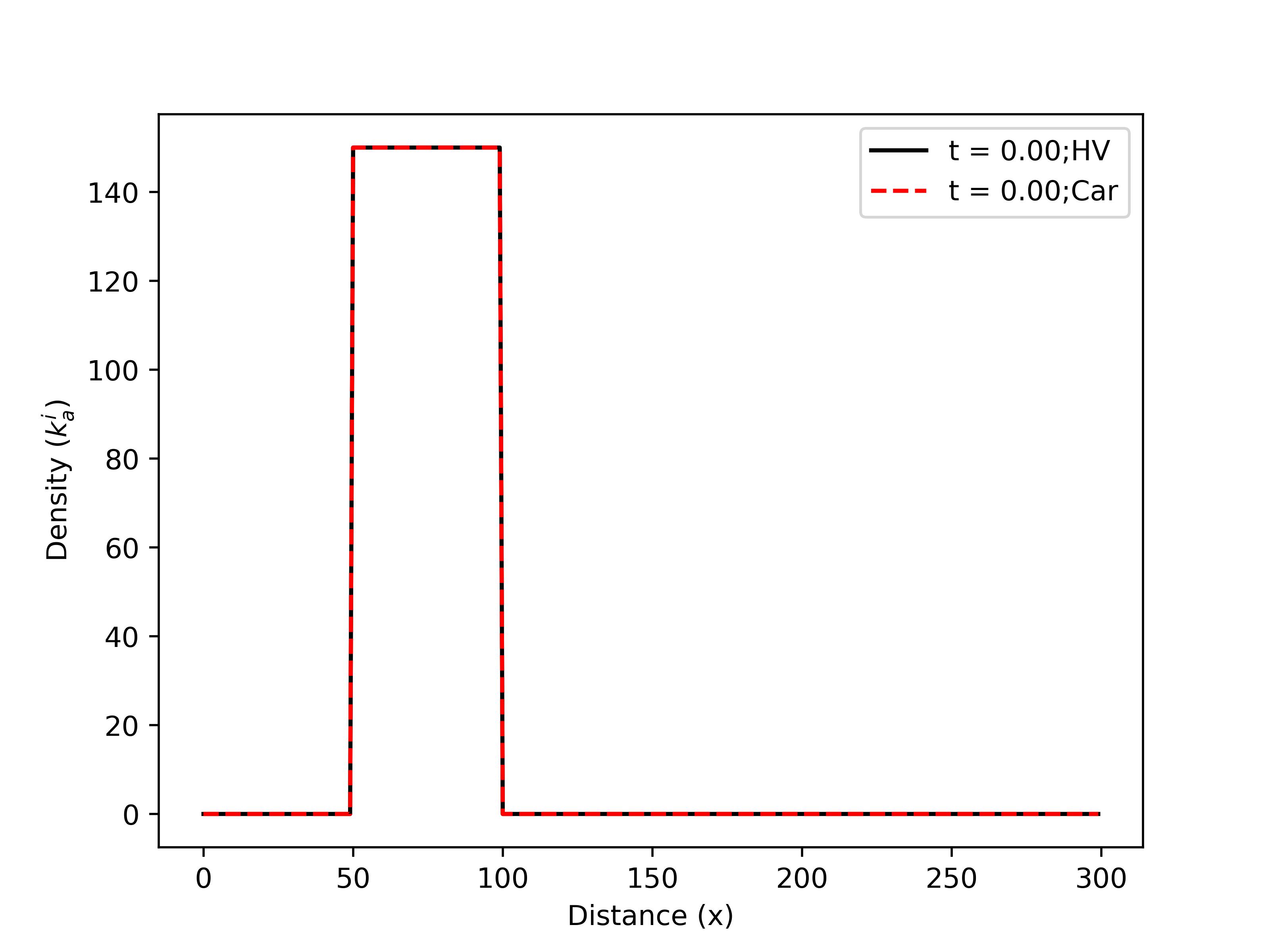}}
\subfigure[$t=80$]{\includegraphics[width=0.22\linewidth]{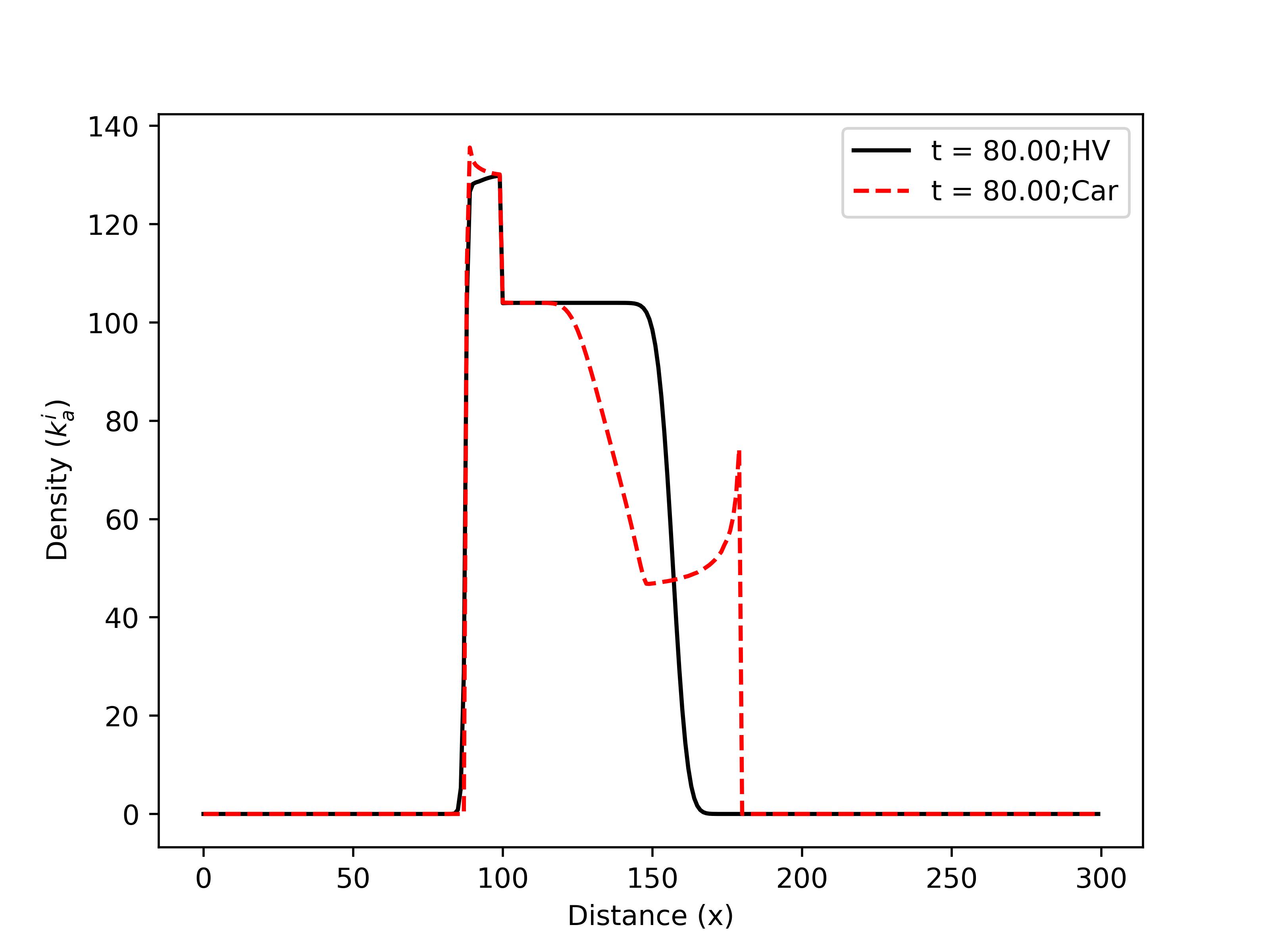}}
\subfigure[$t=160$]{\includegraphics[width=0.22\linewidth]{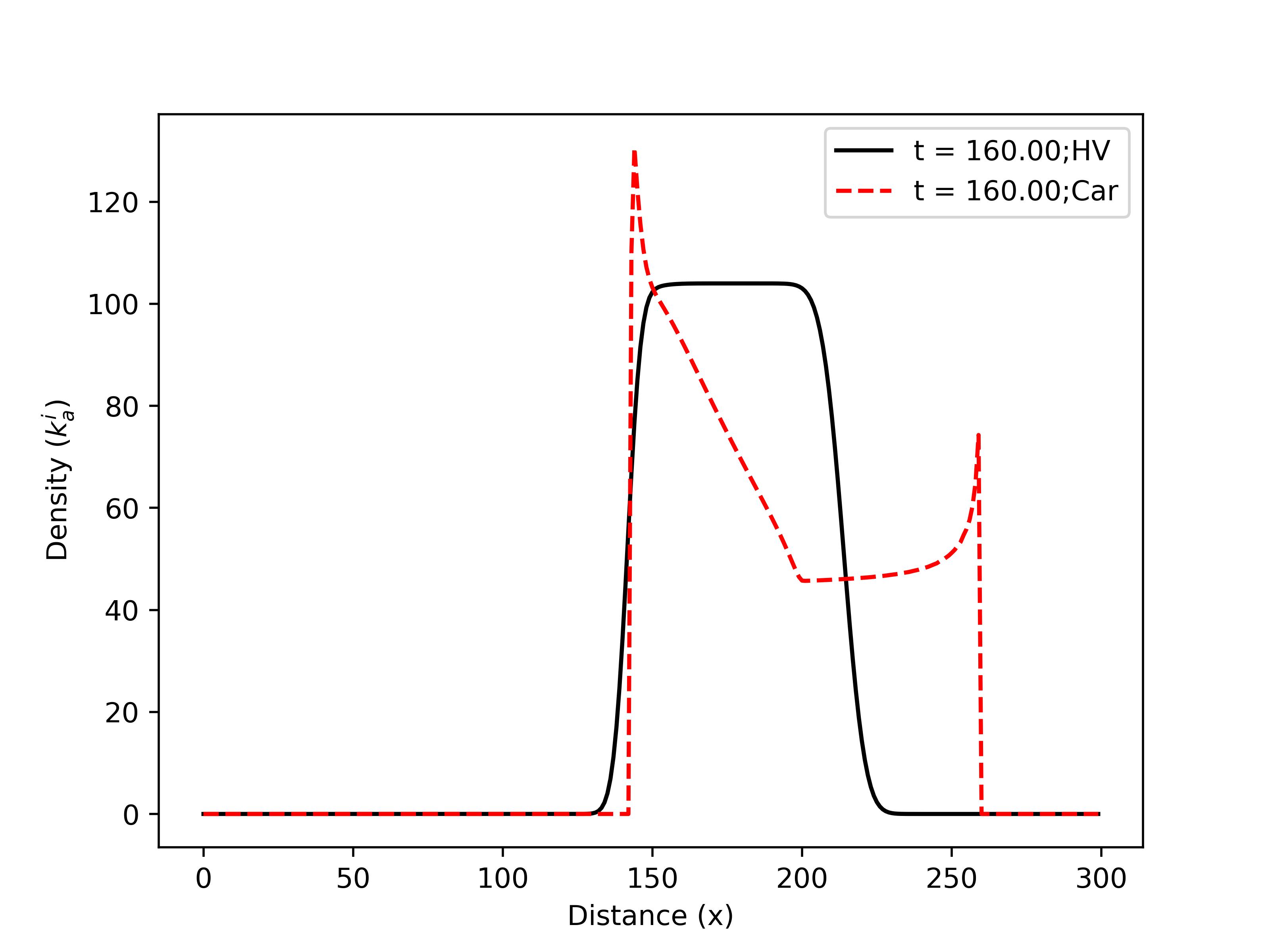}}
\subfigure[Separate Class-specific free flow platoons $(k^{car}_a=150, k^{car}_a=150)$ create congestion after meeting]{\includegraphics[width=0.3\linewidth]{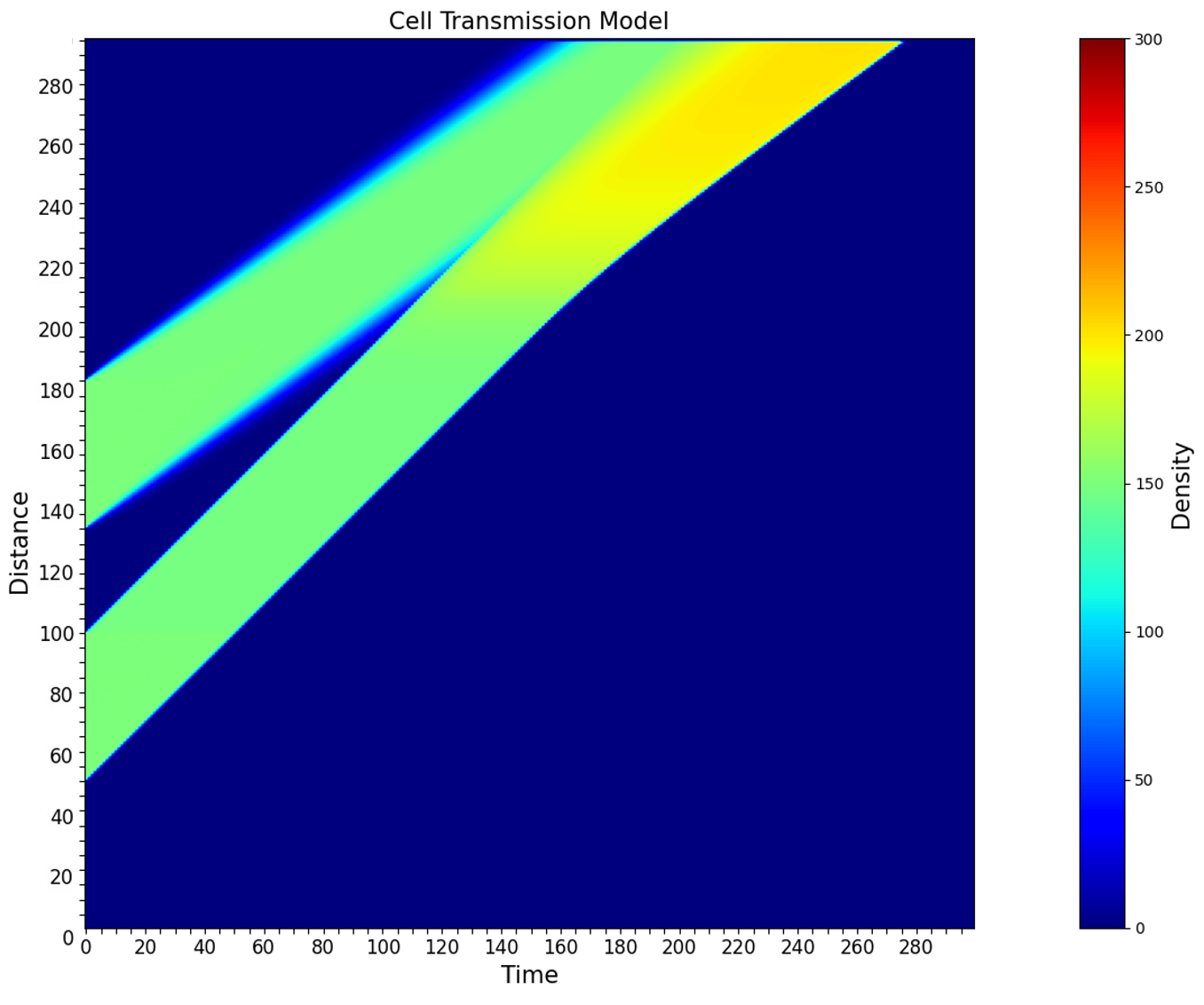}}
\subfigure[$t=0$]{\includegraphics[width=0.22\linewidth]{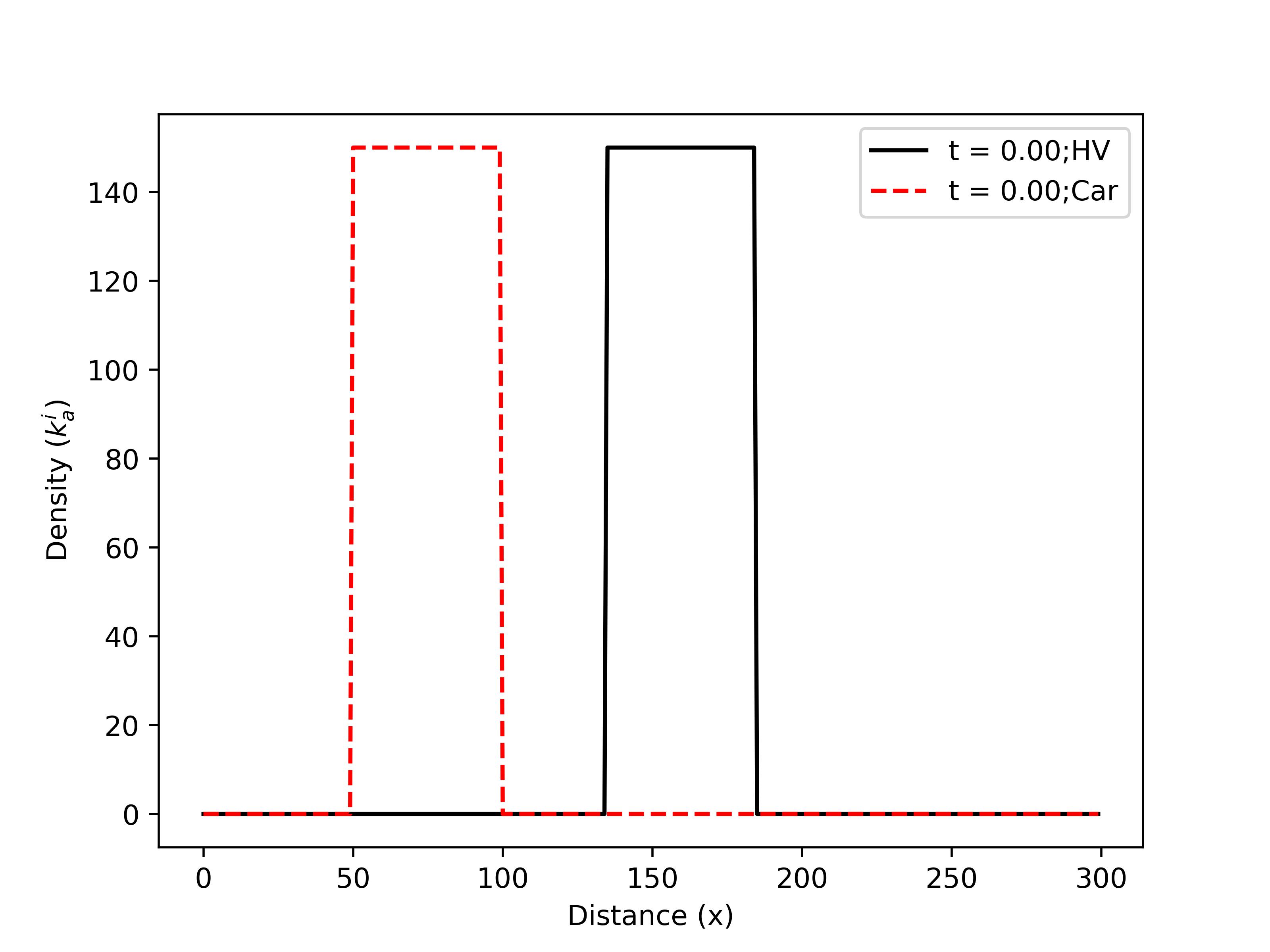}}
\subfigure[$t=80$]{\includegraphics[width=0.22\linewidth]{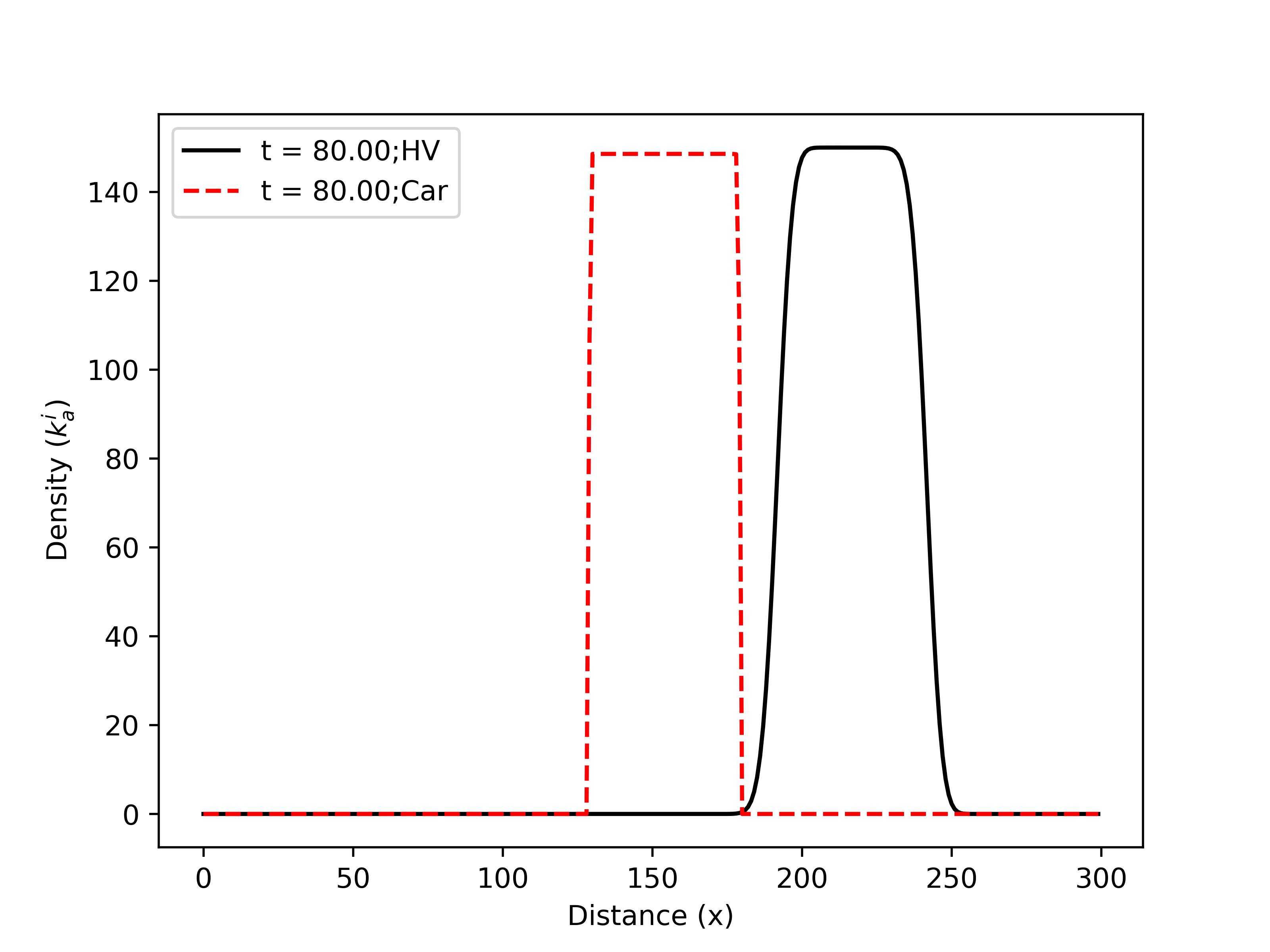}}
\subfigure[$t=140$]{\includegraphics[width=0.22\linewidth]{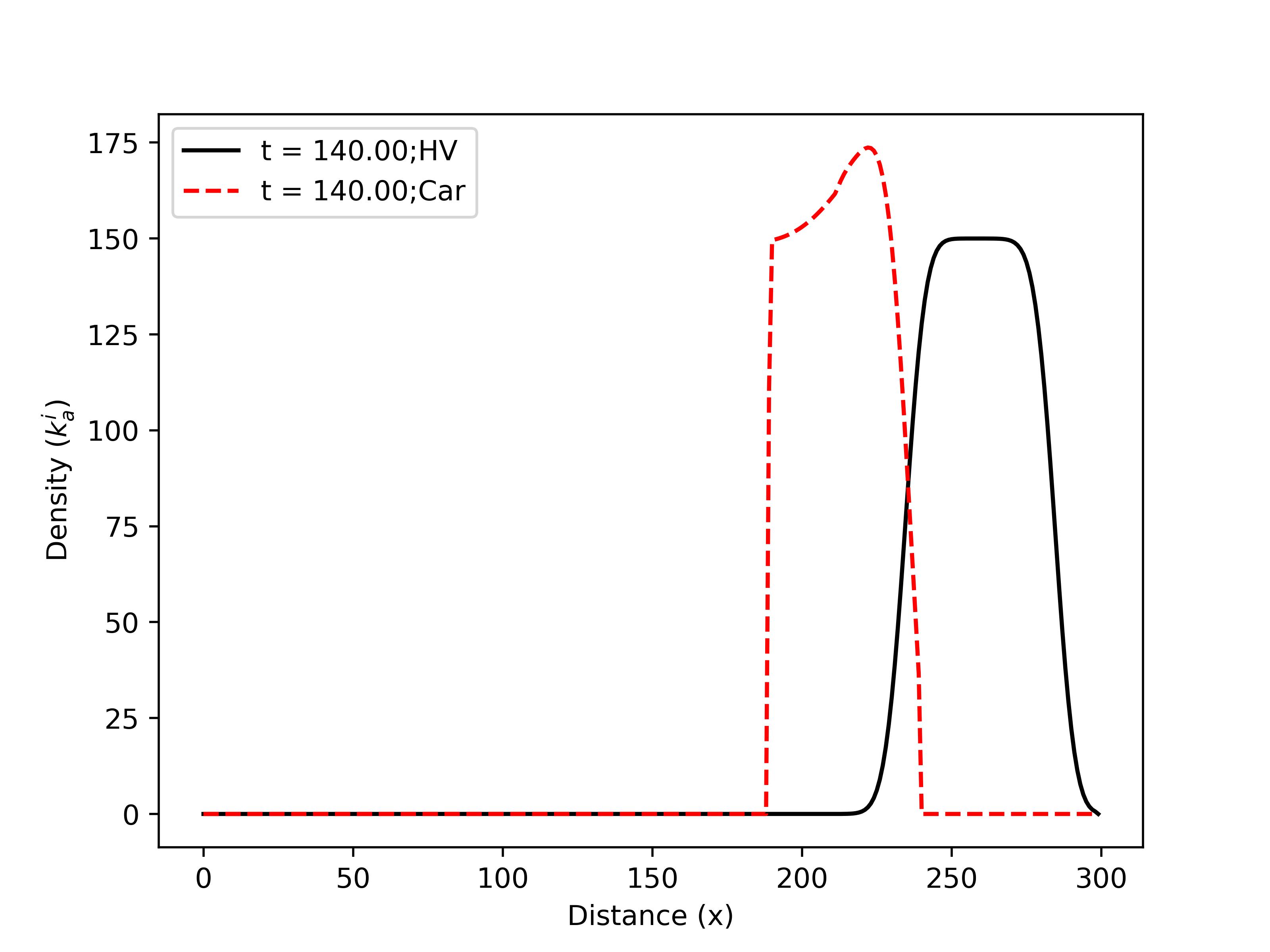}}
\caption{Platoon density propagation and overtaking phenomena of mixed vehicle platoon using proposed m-CTM model with different initial conditions}
\label{fig:mc_ctm_exp}
\end{figure}

This section explores platoon dispersion using the proposed areal continuum model. The initial conditions involve a platoon comprising mixed-class vehicles on a highway, specifically cars and heavy vehicles (HVs). This numerical approach employs the calibrated Chennai multiclass fundamental diagrams and their parameters, as detailed in Table \ref{tab:model_com_smulders}. To illustrate the dispersion of high-speed vehicles from upstream to downstream, we assume that the free-flow speed of cars is 50 km/h, whereas for HVs, it is 45 km/h. The critical density is 200 \(m^2/km-m\) for cars and 250 \(m^2/km-m\) for HVs, resulting in critical speeds of 25 km/h and 23 km/h, respectively. It is noteworthy that the congested wave speed is higher for HVs at 7.7 km/h compared to 5.25 km/h for cars. The platoon dispersion problem is addressed using the proposed multiclass Cell Transmission Model (m-CTM) scheme. We conducted simulations for two distinct initial traffic scenarios to demonstrate platoon dispersion and overtaking maneuvers within a mixed-class platoon consisting of fast and slow vehicles. In Figure \ref{fig:mc_ctm_exp} (a), the visual representation depicts the initial configuration of cars and HVs spanning the spatial extent of $50 m  \leq x \leq 100 m$. The initial densities for both platoons are characterized by a uniform density of 150 $m^2/km-m$, as illustrated in Figure \ref{fig:mc_ctm_exp} (b). As the total initial density exceeded the critical density for both vehicle classes ($k_{a,cr}^{cars}=200,k_{a,cr}^{HVs}=250$), we observed a backward-moving shock, and the traffic density propagated downstream with the critical density value and then dispersed. As the simulation progresses, the plot tracks the temporal and spatial evolution of the vehicles. In Figures \ref{fig:mc_ctm_exp} (c) and \ref{fig:mc_ctm_exp}(d), the temporal evolution of platoon densities is showcased at time instances $t=80 s$ and $t=160 s$, respectively. Notably, at these points in time, the cars have advanced the HVs. The spatial coverage of the car platoon spans the range $140 \leq x \leq 260$, while the HVs occupy the spatial range $130 \leq x \leq 225$ at the 160-second mark. This observation underscores the dispersion of both vehicle classes over time and space, where high-speed vehicles execute overtaking maneuvers to surpass the slow-moving vehicles. The backward shock-wave generated and propagated upstream at $x=100 m$ with an approximate areal density of 130 $m^2/km-m$ for both vehicle classes.
\par We consider an additional initial scenario involving a uniform density of 150 $m^2/km-m$. In this configuration, two distinct platoons were present: a cars platoon spanning the spatial domain of $50 m \leq x \leq 100 m$ and an HVs platoon distributed within the region $135 m \leq x \leq 180 m$. The arrangement is visually represented in Figure \ref{fig:mc_ctm_exp}(f). At the specific time instance of $t=130 s$, the cars and HVs platoons converged, prompting the cars to commence overtaking maneuvers to surpass the slower-moving HVs platoon. This overtaking process is highlighted in Figure \ref{fig:mc_ctm_exp}(g) as the cars integrate themselves at the trailing end of the HV platoon. This integration enables the cars to initiate overtaking actions against the slower HVs, as demonstrated in Figure \ref{fig:mc_ctm_exp}(h). It is important to note that both platoons maintained their respective free flow speeds before their meeting. This led to limited dispersion of their densities up to the point of convergence. However, as the cars joined the tail of the HV platoon, the situation changed. The HVs began to disperse, while the cars encountered congestion. Consequently, the car density escalated to 175  $m^2/km-m$ by the 140-second mark. 


Therefore, Figure \ref{fig:mc_ctm_exp} illustrates the distribution of vehicle density and the shock waves within the platoon over different time intervals along the highway. Over time, it becomes evident that the platoon's leading and trailing ends exhibit a dispersal pattern. This phenomenon can be easily attributed to the presence of drivers with varying behaviors. When opportunities for overtaking arise, vehicles with higher desired free speeds (representing cars) tend to move toward the front of the platoon. Conversely, vehicles with slower desired free speeds (representing HVs) tend to remain toward the tail. The specific shape of the dispersed platoon is influenced by the distribution of desired free speeds employed in the model. Consequently, the presented model and its numerical scheme effectively elucidate the phenomena of platoon dispersion and overtaking maneuvers.

\section{Conclusion}
The new continuum model resembles the LWR model on several accounts but introduces a fresh perspective by considering the conservation of vehicle area instead of the vehicle number.
While a direct comparison with the LWR model wasn't conducted, as the two models represent traffic within different frameworks. The traditional framework uses the number of vehicles, while the proposed method utilizes the projected vehicle area. As mentioned in the introduction, in mixed traffic scenarios, the number of vehicles is less informative without considering vehicle composition. In such cases, using area as a metric is more appropriate. Representing vehicle density or flow as a percentage of the available road area provides a clearer understanding of mixed traffic states compared to using the number of mixed vehicles. Therefore, we believe that the proposed method offers a new perspective on mixed traffic modeling. Although we do not claim that the proposed method is more efficient than the existing model, all existing numerical methods integrate seamlessly with the new model, indicating its compatibility and potential utility in traffic analysis.

\par This shift in focus provides a better representation of traffic dynamics, highlighting the spatial interactions and physical extent of vehicles on the road. The application of vehicle area conservation leads to deriving new fundamental variables that enhance our understanding of mixed traffic flow. These variables are defined based on the conserved vehicle area and offer valuable insights into the macroscopic behavior of traffic.
The newly derived fundamental variables based on conserved vehicle area offer several advantages and implications for traffic flow modeling:
\begin{enumerate}
    \item The new approach acknowledges the spatial aspects of traffic flow, emphasizing the physical interactions and congestion at a macroscopic level. Traditional traffic flow models often assume uniform vehicle sizes, neglecting the impact of heterogeneous vehicle dimensions. The conserved vehicle area-based approach inherently accounts for size variations, leading to a more realistic representation of traffic dynamics.
    \item The new fundamental variables provide a fresh perspective for traffic flow optimization. The area-based variables offer additional metrics to assess traffic performance and identify potential bottlenecks. Traditionally, traffic flow is quantified as the number of vehicles passing through a point per unit of time. With the new approach, traffic flow is defined as the total vehicle area passing through the unit width of the road per unit of time. This parameter provides a spatial measure of traffic flow and incorporates the physical size of vehicles.
    \item Generalized definitions of the proposed variables, following the framework of \citet{edie1963discussion}, enhance confidence in the applicability and generalization of these variables across various traffic studies
    \item  Density, a fundamental measure of traffic concentration, is redefined in terms of the area occupied by vehicles per unit width per unit length of the road rather than the number of vehicles per unit length. Areal density provides a more intuitive representation of traffic density, accounting for the varying vehicle sizes and spatial distribution.
\end{enumerate}
This empirical study examined the correlations among the recently developed areal density, areal flow, and speed using multiple datasets. The findings derived from the empirical data indicate the presence of the Smulders fundamental relationships among the variables. Specifically, these relationships exhibit two distinct regimes: a parabolic pattern during free-flow conditions and a linear pattern during congestion in areal flow-density relationships. These empirical results align with established traffic flow theories and provide substantial evidence supporting the validity and applicability of the newly introduced area-based variables. Additionally, the study investigated the impact of the multi-class modeling approach on traffic flow dynamics along urban inter-city corridors. The findings indicate that the model can explain particular traffic dynamics, including platoon dispersion and freeway traffic behavior. Regarding solution methods for the continuum model, we investigated the applicability of the most popular analytical solution method as the method of characteristics and Godunov's scheme for numerical computation. Since the proposed areal continuum model falls under the class of hyperbolic PDEs, all the available solutions for hyperbolic PDEs seem to apply to the solution. 

Some of the major insights/highlights of the paper are:
\begin{enumerate}
\item A continuum model for heterogeneous traffic was derived based on the principle of vehicle area conservation, and traffic flow parameters based on vehicle area, namely areal flow and areal density, were defined.
\item Equilibrium states for the newly developed parameters were identified using cumulative areal flow and occupancy time curves. A Smulders bi-variate relationship was established between the proposed area-based traffic parameters, and the analysis suggests the stream flow models may better represent mixed traffic than the class-specific representation.
\item The developed model was validated by comparing it with traditional static models using empirical heterogeneous traffic data.
\item A numerical solution to the multi-class kinematic wave model (m-CTM) was proposed to explain various traffic flow phenomena observed in heterogeneous urban traffic.
\item The shock speed and wave speed were described analytically as part of the model's analysis.
\end{enumerate}
It is worth noting that while the empirical analysis provides valuable insights, further research is required to validate and generalize these findings across diverse traffic conditions and geographical contexts. Additionally, future studies could explore the potential applications of the multi-class modeling approach in understanding and managing specific traffic phenomena.
\section*{Acknowledgement}
We would like to express our sincere gratitude to Prof. C. Mallikarjuna (IIT Guwahati), Prof. S Arkatkar (NIT Surat), and Dr. P V Suvin (IIT Hyderabad)  for sharing the data sets. Additionally, we would like to acknowledge the NSFG grant, IIT Madras, International Immersion Exchange Program (IIE), IIT Madras, and Georgia Institute of Technology, Atlanta, for supporting this research.

\makenomenclature
\nomenclature{\(e_i,b_i\)}{Vehicle area comes in}
\nomenclature{\(m,n\)}{Vehicle number comes in}
\nomenclature{\(c_j,d_j\)}{Vehicle area goes out}
\nomenclature{\(p,q\)}{Vehicle number goes out}
\nomenclature{\(\delta\)}{Source or sink vehicle area}
\nomenclature{\(\Delta_x,\Delta_t\)}{Space, and time step}
\nomenclature{\(w\)}{Road width variables}
\nomenclature{\(q_a\)}{Areal flow}
\nomenclature{\(k_a\)}{Areal density}
\nomenclature{\(v_a\)}{Areal speed}
\nomenclature{\(N\)}{Number of vehicles}
\nomenclature{\(W\)}{Road width}
\nomenclature{\(X\)}{Road length}
\nomenclature{\(T\)}{Time }
\nomenclature{\(O_c\)}{Occupancy}
\nomenclature{\(ao\)}{Area occupancy}
\nomenclature{\(L_i\)}{Length of the $i^{th}$ vehicle}
\nomenclature{\(B_i\)}{Width of the $i^{th}$ vehicle}
\nomenclature{\(d\)}{Length of the detector}
\nomenclature{\(TW\)}{Two-wheelers}
\nomenclature{\(ThW\)}{Three-wheelers}
\nomenclature{\(LCV\)}{Light-Commercial-Vehicles}
\nomenclature{\(HV\)}{Heavy Vehicles}
\nomenclature{\(A(X,t)\)}{Cumulative vehicle area in a location $X$ at time $t$}
\nomenclature{\(a_0,b_0\)}{average areal flow and areal occupancy}
\nomenclature{\(k_{a,jam}\)}{Jam areal density}
\nomenclature{\(k_{a,crit}\)}{Critical areal density}
\nomenclature{\(v_{crit}\)}{Critical speed}
\nomenclature{\(v_{max}\)}{Maximum speed}
\nomenclature{\(v_{f}\)}{Free flow speed}
\nomenclature{\(\omega_{a}\)}{areal wave speed at congestion}
\nomenclature{\(c(k_a)\)}{wave speed as a function of areal density $k_a$}
\nomenclature{\(s\)}{slope of the wave speed line}
\nomenclature{\(\mu(k_a)\)}{Supply function}
\nomenclature{\(\lambda(k_a)\)}{Demand function}
\nomenclature{\(\Phi\)}{Average flow in CTM model}
\printnomenclature

\appendix
\section{APPENDIX}
\begin{table}[H]
\caption{Comparison of model parameters for two-wheelers with the model fitness score for all data}
\label{tab:model_com_tw}
\begin{tabular}{@{}llllllllll@{}}
\toprule
\multicolumn{6}{c}{Two-wheelers (TWs)}                                     & \multicolumn{2}{c}{$k_a-v^i_a$}       & \multicolumn{2}{c}{$q^i_a-k_a$}      \\ \midrule
Location & Model               & $v^i_{f}$ & $v^i_{cr}$ & $k^i_{a,cr}$ & $\omega^i_a$   & $R^2$         & $RMSE$          & $R^2$         & $RMSE$         \\ \midrule
         & Greenshields               & 37   &   -   &    -   &    - & 0.54          & 8.60           & 0.07          & 2236         \\
         & Greenberg                  &   -   & 14.3 &  -     &    - & 0.76          & 7.2           & 0.62          & 950          \\
Chennai  & Underwood                  & 51   &   -   & 284   &   -  & 0.83          & 5.70           & 0.64          & 927          \\
         & Del Castillo               & 41.7 & -     &   -    & 6.3 & 0.82          & 5.43          & 0.62          & 971          \\
         & Daganzo                    & 40   &   -   & 150   & 6.6 & 0.78          & 6.20           & 0.54          & 1057         \\ 
         & Smulders                    & 49.5   &   29   & 170   & 5.94 & 0.84          & 5.32           & 0.62         & 973        \\\hline
         & Greenshields               & 36.7 &   -   &    -   &  -   & 0.64          & 8.52          & 0.08          & 2046         \\
         & Greenberg                  &    -  & 12   &   -    &  -   & 0.82          & 6.12          & 0.62          & 720          \\
Surat    & Underwood                  & 48   &    -  & 248   &  -   & 0.88          & 4.70           & 0.72          & 614          \\
         & Del Castillo               & 40.5 &    -  &    -   & 4.6 & 0.89          & 5.10           & 0.62          & 725          \\
         & Daganzo                    & 39   &   -   & 155   & 4   & 0.79          & 6.42          & 0.62          & 737          \\ 
         & Smulders                    & 48.40  &   19.5   & 193   & 4.66 & 0.90          & 4.53           & 0.74         & 603        \\\hline
         & Greenshields               & 37   &   -   &     -  &   -  & 0.72          & 7.68          & 0.09          & 2295         \\
         & Greenberg                  &    -  & 13.2 &  -     &   -  & 0.83          & 5.90           & 0.63          & 851          \\
Guwahati & Underwood                  & 52   &    -  & 290   &    - & 0.89          & 4.50           & 0.69          & 817          \\
         & Del Castillo               & 40.5 &   -   &     -  & 7.5 & 0.89          & 4.10           & 0.79          & 661          \\
         & Daganzo                    & 38.5 &    -  & 155   & 6.8 & 0.85          & 4.95          & 0.72          & 768          \\ 
         & Smulders                    & 48.80   &   29.40   & 170   & 6.02 & 0.90          & 3.93          & 0.80        & 647        \\\bottomrule
\end{tabular}
\end{table}

\begin{table}[H]
\caption{Comparison of model parameters for cars with the model fitness score for all data}
\label{tab:model_com_car}
\begin{tabular}{@{}llllllllll@{}}
\toprule
\multicolumn{6}{c}{Cars}                                     & \multicolumn{2}{c}{$k_a-v^i_a$}       & \multicolumn{2}{c}{$q^i_a-k_a$}      \\ \midrule
Location & Model               & $v^i_{f}$ & $v^i_{cr}$ & $k^i_{a,cr}$ & $\omega^i_a$   & $R^2$         & $RMSE$          & $R^2$         & $RMSE$         \\ \midrule
         & Greenshields               & 35   &   -    &    -   &    - & 0.56          & 8.62          & 0.05          & 2240         \\
         & Greenberg                  &    -  & 13.4  &    -   &    - & 0.79          & 6.90          & 0.63          & 955          \\
Chennai  & Underwood                  & 49   &    -   & 265   &   -  & 0.83          & 5.42          & 0.65          & 928          \\
         & Del Castillo               & 39   &    -   &   -    & 6.2 & 0.83          & 5.42          & 0.63          & 978          \\
         & Daganzo                    & 37   &       & 145   & 6.2 & 0.79          & 6.07          & 0.56          & 1049         \\ 
         & Smulders                    & 49   &   25   & 200   & 5.25 & 0.82          & 6.50           & 0.61         & 980       \\\hline
         & Greenshields               & 35.2 &   -    &    -   &   -  & 0.65          & 8.50           & 0.07          & 2055         \\
         & Greenberg                  &    -  & 11.65 &   -    &    - & 0.83          & 5.85          & 0.63          & 710          \\
Surat    & Underwood                  & 46   &    -   & 237   &   -  & 0.89          & 4.60           & 0.73          & 615          \\
         & Del Castillo               & 39.4 &   -    &    -   & 4.2 & 0.88          & 4.90           & 0.64          & 722          \\
         & Daganzo                    & 42   &    -   & 160   & 4.8 & 0.78          & 6.48          & 0.61          & 742          \\ 
         & Smulders                    & 48.80   &   19.70   & 195   & 4.77 & 0.89          & 4.80           & 0.73         & 612       \\\hline
         & Greenshields               & 36.4 &    -   &     -  &   -  & 0.71          & 7.64          & 0.08          & 2289         \\
         & Greenberg                  &    -  & 13.1  &    -   &  -   & 0.84          & 5.20           & 0.67          & 843          \\
Guwahati & Underwood                  & 49.5 &    -   & 268   &    - & 0.90           & 4.00             & 0.70           & 810          \\
         & Del Castillo               & 37.8 &  -     &     -  & 6.4 & 0.90           & 4.07          & 0.80           & 659          \\
         & Daganzo                    & 35.5 &     -  & 160   & 6.2 & 0.86          & 4.89          & 0.73          & 765          \\ 
         & Smulders                    & 48.80  &   25.60  & 195   & 6.20 & 0.91         & 3.92           & 0.81         & 643       \\\bottomrule
\end{tabular}
\end{table}

\begin{table}[H]
\caption{Comparison of model parameters for heavy vehicles with the model fitness score for all data}
\label{tab:model_com_hvs}
\begin{tabular}{@{}llllllllll@{}}
\toprule
\multicolumn{6}{c}{Heavy vehicles (HVs)}                                     & \multicolumn{2}{c}{$k_a-v^i_a$}       & \multicolumn{2}{c}{$q^i_a-k_a$}      \\ \midrule
Location & Model               & $v^i_{f}$ & $v^i_{cr}$ & $k^i_{a,cr}$ & $\omega^i_a$   & $R^2$         & $RMSE$          & $R^2$         & $RMSE$         \\ \midrule
        & Greenshields               & 34   &    -  &    -   &    -  & 0.59          & 8.58          & 0.09          & 2252         \\
         & Greenberg                  &    -  & 13.6 &   -    &   -   & 0.8           & 6.00            & 0.64          & 947          \\
Chennai  & Underwood                  & 51.4 &   -   & 260   &   -   & 0.82          & 5.60           & 0.65          & 930          \\
         & Del Castillo               & 40.7 &   -   &    -   & 6.4  & 0.82          & 5.44          & 0.63          & 970          \\
         & Daganzo                    & 40   &    -  & 135   & 6.3  & 0.79          & 6.11          & 0.55          & 1055         \\ 
         & Smulders                    & 49   &   23  & 250   & 7.67 & 0.84          & 5.40           & 0.62         & 975       \\\hline
         & Greenshields               & 35   &    -  &     -  &   -   & 0.66          & 8.53          & 0.06          & 2015         \\
         & Greenberg                  &  -   & 13.2 &    -   &    -  & 0.83          & 5.87          & 61.00           & 715          \\
Surat    & Underwood                  & 46   &    -  & 235   &   -   & 0.87          & 4.80           & 0.71          & 620          \\
         & Del Castillo               & 43   &   -   &   -    & 4    & 0.87          & 5.00             & 0.63          & 721          \\
         & Daganzo                    & 38   &      & 140   & 5.8  & 0.77          & 6.50           & 0.62          & 740          \\ 
         & Smulders                    & 48.2   &   20   & 210   & 5.32 & 0.89         & 4.60           & 0.72         & 604       \\\hline
         & Greenshields               & 37   &   -   &   -    &   -   & 0.73          & 7.70           & 0.07          & 2275         \\
         & Greenberg                  &     - & 15   &    -   &   -   & 0.81          & 5.30           & 0.66          & 845          \\
Guwahati & Underwood                  & 52   &    -  & 264   &    -  & 0.88          & 4.20           & 0.68          & 825          \\
         & Del Castillo               & 39.5 &    -  &   -    & 6.54 & 0.90           & 4.11          & 0.80           & 660          \\
         & Daganzo                    & 37.5 &   -   & 140   & 5.8  & 0.84          & 4.90           & 0.72          & 765          \\ 
         & Smulders                    & 48.65   &   30   & 200   & 7.50 & 0.89          & 3.96           & 0.79         & 645       \\\bottomrule
\end{tabular}
\end{table}


\newpage
\bibliographystyle{elsarticle-num-names} 
\bibliography{cas-refs}

\end{document}